\documentclass{amsart}

\usepackage[english]{babel}
\usepackage{graphicx}
\usepackage[hidelinks]{hyperref} 
\usepackage{amssymb}
\usepackage{amsmath}
\usepackage[foot]{amsaddr}
\usepackage{amsfonts}
\usepackage{xcolor}
\usepackage{mathtools}
\usepackage{algpseudocode}
\usepackage{algorithm}
\usepackage{comment}
\usepackage{subcaption}
\usepackage{caption}
\usepackage[top=1in, bottom=1.25in, left=1.25in, right=1.25in]{geometry}
\usepackage{float}
\usepackage{tikz}
\newcommand{\norm}[1]{\lVert #1 \rVert}

\newcommand{\bmu}{{{\mu}}}
\newcommand{\V}{{\mathbb{V}}}
\newcommand{\Q}{{\mathbb{Q}}}
\newcommand{\Cal}[1]{{\mathcal{#1}}}
\newcommand{\mbb}[1]{{\mathbb{#1}}}

\DeclareMathOperator*{\argmax}{arg\,max}

\DeclareMathOperator*{\argmin}{arg\,min}

\usepackage{todonotes} 
\newcommand{\B}[1]{{ \color{black}{#1}}}
\newcommand{\A}[1]{{ \color{black}{#1}}}

\makeatletter
\algnewcommand{\LineComment}[1]{\Statex \hskip\ALG@thistlm {\color{white}\textbf{Input:}} #1}
\makeatother

\makeatletter
\algnewcommand{\LineCommentO}[1]{\Statex \hskip\ALG@thistlm {\color{white}\textbf{Output:}} #1}
\makeatother
{}

\newtheorem{remark}{Remark}

\algrenewcommand\algorithmicrequire{\textbf{Input:}}
\algrenewcommand\algorithmicensure{\textbf{Output:}}
\usepackage[toc,page]{appendix}

\begin{document}
\title[Deflated-Greedy and adaptivity for certified bifurcations]{Deflation-based certified greedy algorithm and adaptivity for bifurcating nonlinear PDEs}

\author{Federico Pichi$^{1,2, *}$ and Maria Strazzullo$^{3, *}$}
\address{$^1$ Chair of Computational Mathematics and Simulation Science, École Polytechnique Fédérale de Lausanne, 1015 Lausanne, Switzerland}
\address{$^2$ mathLab, Mathematics Area, SISSA, via Bonomea 265, I-34136 Trieste, Italy}
\address{$^3$ Politecnico di Torino, Department of Mathematical Sciences ``Giuseppe Luigi Lagrange'', Corso Duca degli Abruzzi, 24, 10129, Torino, Italy}
\address{$^*$ INdAM-GNCS group member}

\begin{abstract}
This work deals with tailored reduced order models for bifurcating nonlinear parametric partial differential equations, where multiple coexisting solutions arise for a given parametric instance. Approaches based on proper orthogonal decomposition have been widely investigated in the literature. Still, they usually rely on some \emph{a priori} information about the bifurcating model and lack any error estimation. On the other hand, standard certified reduced basis techniques fail to represent correctly the branching behavior, given that the error estimator is no longer reliable. The main goal of the contribution is to overcome these limitations by introducing two novel iterative algorithms, namely: (i) the adaptive-greedy, detecting the bifurcation point starting from scarce information over the parametric space, and (ii) the deflated-greedy, certifying multiple coexisting branches simultaneously. The former approach takes advantage of the non-smoothness of the reduced manifold to detect the bifurcation, while the latter exploits the deflation and continuation methods to enrich the reduced space with the bifurcating solutions.
We test the two strategies for the Coanda effect held by the Navier-Stokes equations in a sudden-expansion channel. The accuracy and error certification are compared with standard greedy and proper orthogonal decomposition.

\end{abstract}

\maketitle
\tableofcontents

\section{Introduction and motivation}

In many scientific and industrial applications, nonlinear parametric partial differential equations (PDEs) are employed as models to reliably describe
complex physical phenomena. 
These models play a ubiquitous role in many fields, in contrast with linear equations that provide a simplified setting, possibly discarding fundamental information regarding the model's behavior, e.g., the coexistence of multiple solutions in bifurcating scenarios \cite{ambrosetti1995primer,CalozNumericalAnalysisNonlinear1997,seydel2009practical,KielhoferBifurcationTheoryIntroduction2012,KuznetsovElementsAppliedBifurcation2023}.
The investigation of such problems poses different challenges already at the theoretical level: from the detection and analysis of the bifurcating branches and their features to the study of the stability properties of the solutions. Indeed, these bifurcating models are characterized by drastic changes in the system's state as a consequence of a slight variation of the parameter. This non-differentiable evolution of the solution with respect to (w.r.t.)\ the parametric nature of the problem gives rise to the non-uniqueness of the solution of the PDE. Several physical phenomena are described by well-known models that suffer from this ill-posedness, from beams' buckling in computational mechanics \cite{bigoni2012nonlinear,PichiReducedOrderModels2023,NiroomandiModelOrderReduction2010,PichiReducedBasisApproaches2019} to symmetry-breaking flows in fluid dynamics \cite{PittonComputationalReductionStrategies2017,QuainiSymmetryBreakingPreliminary2016,PichiDrivingBifurcatingParametrized2022a,PichiArtificialNeuralNetwork2023,DijkstraBifurcationAnalysisFluid2023a,BattagliaBifurcationLowReynolds1997}, passing through particles' configuration in quantum systems \cite{CharalampidisComputingStationarySolutions2018,CharalampidisBifurcationAnalysisStationary2020,PichiReducedOrderModeling2020} and Rayleigh–Bénard convection problems \cite{BoulleBifurcationAnalysisTwodimensional2022,LaakmannBifurcationAnalysisTwodimensional2024}, possibly affected by stochastic behaviors \cite{VenturiStochasticBifurcationAnalysis2010,GonnellaStochasticPerturbationApproach2024a,KuehnUncertaintyQuantificationAnalysis2024,KuehnUncertaintyQuantificationBifurcations2021,JinNovelStochasticBifurcation2022}.
The computational analysis of bifurcating PDEs comes with many difficulties. Indeed, to recover the so-called bifurcation diagram, i.e., a comprehensive qualitative analysis of the branching behavior, might require (i) a fine mesh resolution to capture the effect of nonlinear terms, (ii) many-query parametric evaluations and (iii) tailored methodologies to extract the bifurcating information. 
In this context, standard discretization techniques, i.e., the high-fidelity (HF) approaches based, e.g., on the Finite Element (FE) methods, might lead to unaffordable
computational costs, especially when considering design and optimization analysis settings. For this reason, there is a growing interest in developing Reduced Order Models (ROMs) \cite{hesthavenCertifiedReducedBasis2015,QuarteroniReducedBasisMethods2016} capable of dealing with bifurcating PDEs depending on parameters that might change the physics or the geometry of the system at hand \cite{BruntonDiscoveringGoverningEquations2016,DengLoworderModelSuccessive2020a,HerreroRBReducedBasis2013,HessLocalizedReducedorderModeling2019,pichi_phd}.
The purpose of ROMs is to build a surrogate low-dimensional model that approximates the HF solution efficiently and reliably.

In the literature, different approaches have been investigated to build such reduced models, ranging from projection-based techniques to data-driven approaches \cite{PittonApplicationReducedBasis2017,PichiReducedBasisApproaches2019,PintoreEfficientComputationBifurcation2021,PichiArtificialNeuralNetwork2023,PichiGraphConvolutionalAutoencoder2024,HessLocalizedReducedorderModeling2019,SilvesterMachineLearningHydrodynamic2024,FabianiNumericalSolutionBifurcation2021,OlshanskiiApproximatingBranchSolutions2025}. When choosing among these two categories for a specific application, one should consider the trade-off between the reliability and robustness of the former and the efficiency of the latter. Here, we focus on projection-based approaches, providing a mathematically consistent technique for a certified recovery of the bifurcating phenomena. 

Projection-based techniques are usually based on the selection, computation, and successive exploitation of a dataset constituted by HF solutions, called \emph{snapshots}. The two main approaches employed to build the reduced space are (i) the Proper Orthogonal Decomposition (POD) \cite{hesthavenCertifiedReducedBasis2015}, which performs a data compression extracting the most meaningful modes representing the HF solution, and (ii) the greedy algorithm, which adaptively enriches the basis of the reduced problem by exploiting an error estimator \cite{buffa2012priori,QuarteroniReducedBasisMethods2016}. Historically, the POD approach has been preferred in the bifurcating setting, either compressing the whole bifurcating scenario (global POD) or each specific branch separately (branch-wise POD) \cite{pichi_phd}. 

These techniques can only recover the \emph{sampled} bifurcating states. Namely, the HF bifurcating snapshots are needed to build representative modes and, if no coexisting solutions are detected, the bifurcation information is lost and the POD fails in representing the bifurcation diagram. Moreover, this approach does not guarantee an error certification, in general and even more in this specific case, meaning that the approximation accuracy can only be inferred by the decay of the singular values assuming a complete knowledge of the phenomenon.

On the other hand, while the vanilla-greedy approach allows for error certification in non-bifurcating cases, the non-uniqueness of the solution for a given parameter value prevents its reliability. 
In this work, we overcome these issues, proposing two novel algorithms to deal with bifurcating systems. More specifically, we develop tailored reduced basis (RB) greedy strategies for nonlinear PDEs enhanced by the error certification of multiple branches with no prior knowledge of the bifurcating nature of the system.
The main novelties of this contribution are:
\begin{itemize}
    \item[$\circ$] the adaptive-greedy strategy, which iteratively samples parameters in a neighborhood of the critical ones at which the bifurcation occurs, even starting from a coarse discretization of the parametric space, automatically 
    detecting the bifurcations points;
    \item[$\circ$] the deflated-greedy strategy, an algorithm that searches for multiple bifurcating snapshots via the deflation method, enriches the reduced space with the bifurcating features of PDE via a reduced deflation method, computes multiple error estimators related to different solution behaviors, and certifies them simultaneously.
\end{itemize} 
The algorithms are tested on a well-known bifurcating problem in fluid dynamics, the so-called Coanda effect, consisting of a sudden-expansion channel flow held by the Navier-Stokes equations featuring a supercritical pitchfork bifurcation \cite{tritton2012physical,QuainiSymmetryBreakingPreliminary2016}. The main strength of the novel approaches is their reliability. Indeed, even being unaware of the possible bifurcating nature of the PDE under investigation, we can certify the error w.r.t.\ the different branches. Moreover, the adaptive-greedy helps in reducing the number of training parameters to initialize the greedy approach. Indeed, considering a large parametric set usually leads to a waste of computational time, by exploring the uniqueness basin and only sampling a few points revealing the bifurcating nature of the system.
To the best of our knowledge, these are the first greedy-based strategies developed to investigate bifurcating PDEs.

The paper is outlined as follows: in Section \ref{sec:ROMnonlinear} we describe the model setting at the continuous and discrete level focusing on certified ROMs under the hypothesis of the uniqueness of the solution. Section \ref{sec:sota} introduces bifurcating nonlinear parametric PDEs, and how they have been classically tackled both at the HF and ROM levels, underlying the limitations of ROM techniques in such context. In Sections \ref{sec:adaptivegreedy} and \ref{sec:defgreedy}, we present the adaptive-greedy and the deflated-greedy approaches, respectively.
The new algorithms are compared with greedy and POD approaches in terms of accuracy w.r.t.\ the HF model on Navier-Stokes equations in a sudden-expansion channel geometry in Section \ref{sec:results}. Conclusions follow in Section \ref{sec:conc}. Additional results about nonlinear error estimation 
can be found in 
\ref{app:estimator}.

\section{Certified ROMs for nonlinear parametric PDEs}
\label{sec:ROMnonlinear}
In this section, we recall the main ingredients for obtaining certified ROMs when dealing with parametric nonlinear PDEs.
These techniques build a surrogate model, which is much cheaper to evaluate than the original HF one, exploiting the projection on a lower dimensional space, and the parametric structure of the model itself.
Thus, ROMs are designed to perform experimental analysis in a faster, but still accurate, way. Applying ROMs to nonlinear PDEs might be particularly challenging even in the uniqueness setting, especially for (i) the availability and reliability of \emph{a-posteriori} error estimators and (ii) the non-affine nature of the system. Indeed, the \textit{affine decomposition assumption} guarantees the ROM \emph{offline-online} paradigm. That is, the reduced space is built \emph{offline} using possibly costly procedures, while \emph{online} it can be efficiently queried for many instances of the parameter value. This ideal scenario is not met by nonlinear problems, where the affine decomposition has to be recovered via hyper-reduction techniques \cite{barrault04:_empir_inter_method,ChaturantabutNonlinearModelReduction2010}. Here, we only focus on the role of the error certification during the offline stage, and we do not tackle the online efficiency issue through hyper-reduction. The latter, involving an approximation of the nonlinear term responsible for the bifurcation, entails a new set of challenges that are beyond the scope of the contribution.

\subsection{Problem formulation: continuous and discrete settings}
Let us consider the parameter space $\mathcal P \subset \mathbb R^P$, for a natural number $P\geq 1$, and an open and bounded regular domain $\Omega \subset \mathbb R^d$, where $d$ denotes its spatial dimension. Let $G: \mathbb U \times \mathcal P \rightarrow \mathbb U^*$, with $\mathbb U$ a suitable Hilbert space, be a parametric nonlinear operator defining a PDE. 

The parametric problem is: given $\mu \in \mathcal P$, find $u = u(\mu) \in \mathbb U$ such that
\begin{equation}
    \label{eq:strong_pde}
        G(u; \mu) = 0 \quad \text{in } \mathbb U^*.
\end{equation}
The weak formulation associated to \eqref{eq:strong_pde} reads: for $\mu \in \mathcal P$, find $u \in \mathbb U$ such that
\begin{equation}
    \label{eq:weak_pde}
        g(u,q; \mu) = 
            \left \langle G(u; \mu), q \right \rangle = 0 \quad \text{for all } q \in \mathbb U,
\end{equation}
where $g: \mathbb U \times \mathbb U \rightarrow \mathbb R$ is the variational nonlinear form, and $\left \langle \cdot, \cdot \right \rangle$ the duality pairing between $\mathbb U^*$ and $\mathbb U$. 

In the discrete setting\footnote{The proposed procedure does not depend on the chosen high-fidelity discretization.}, the problem is formulated as: given $\mu \in \mathcal P$, find $u_h = u_h(\mu) \in \mathbb{U}^{N_h}$, with $\mathbb{U}^{N_h} \subset \mathbb{U}$ a finite dimensional subspace of dimension $N_h$ spanned by a set of basis functions $\{\varphi_i\}_{i=1}^{N_h}$, such that
\begin{equation}
    \label{eq:res_vec}
        \mathsf G(\mathsf u_h; \mu) = 0,
\end{equation}
where $\mathsf G(\mathsf u_h; \mu)$ is the \emph{residual vector} defined as 
$
\mathsf G(\mathsf u_h; \mu)_i = g(u_h, \varphi_i; \mu), 
$
for $i=1, \dots, N_h$, and $\mathsf u_h \in \mathbb R^{N_h}$ is the vector of the unknown coefficients. The nonlinear system \eqref{eq:res_vec} can be solved using the Newton method (see Algorithm \ref{alg:newton}), which $k$-th iteration reads: for a given $\mu \in \mathcal P$ and initial guess $\mathsf u_h^{(0)} \in \mathbb R^{N_h}$, find $\delta \mathsf u_h \in \mathbb R^{N_h}$ such that 
\begin{equation}
\label{eq:Newton_method_k}
\mathsf {Jac}^G(\mathsf u_h^{(k)}; \mu) \delta \mathsf u_h = - \mathsf G(\mathsf u_h^{(k)}; \mu), \quad \text{ and } \quad  \mathsf u_h^{(k+1)} = \mathsf u_h^{(k)} + \delta \mathsf u_h,  
\end{equation}
where $\mathsf {Jac}^G(\mathsf u_h^{(k)}; \mu)\in \mathbb R^{N_h \times N_h}$ is the \emph{Jacobian matrix} of the nonlinear system. 
The procedure is performed until a stopping criterion is reached, i.e.\ the norm of the residual or the norm of the increment are below a given tolerance ($\mathsf {tol}$). 
For the sake of clarity, in this section, we assume that the nonlinear problem is well-posed and that the Jacobian matrix is not singular. It is well-known that such a condition corresponds to the positivity of the \textit{discrete inf-sup} constant $\beta^h (\mu)$ in a neighborhood of $\overline{\mathsf u} \in \mathbb R^{N_h}$ defined as

\begin{equation}
\label{eq:inf_sup_disc}
\beta^h (\mu) = \adjustlimits \inf_{\mathsf{u} \neq 0} \sup_{\mathsf{q} \neq 0} \frac{\mathsf{q}^T \mathsf{Jac}^G(\overline{\mathsf u}; \mu)\ \mathsf{u}}{\norm{\mathsf{u}}_{\mbb U^{N_h}}\norm{\mathsf{q}}_{\mbb U^{{N}_h}}} \qquad \forall \, \mu \in \mathcal P,
\end{equation}
for $\mathsf u, \mathsf q \in \mathbb R^{N_h}$. 
This is a common assumption in the literature, but it is violated in many realistic scenarios governed by complex models, originating interesting phenomena.
Indeed, as we will discuss in Section \ref{sec:sota}, for a given parameter $\mu$, a nonlinear PDE might lose the uniqueness of the solution, admitting the coexistence of qualitatively different states.

\begin{algorithm}
\caption{The $ \mathsf {Newton}$ algorithm}\label{alg:newton}
\begin{algorithmic}[1]
\Statex{\textbf{Input:} residual $\mathsf G$, tolerance $\mathsf {tol}$, initial guess $\mathsf u^{(0)}$, parameter $\mu$}
\Statex{\textbf{Output:} solution of the system $\mathsf u$}
\Statex
\While{$\norm{\mathsf G(\mathsf u_h^{(k)}; \mu)}_{\mathbb U^{N_h}} \geq \mathsf{tol}$}
\State{$\mathsf {Jac}^G(\mathsf u_h^{(k)}; \mu) \delta \mathsf u_h = - \mathsf G(\mathsf u_h^{(k)}; \mu)$} \Comment{Solve linearized system}
\State{$\mathsf u_h^{(k+1)} = \mathsf u_h^{(k)} + \delta \mathsf u_h$} \Comment{Update solution}
\State{$k = k + 1$}
\EndWhile
\end{algorithmic}
\end{algorithm}

\subsection{The greedy algorithm}
\label{sec:greedy}
The objective of this section is to discuss strategies to build a low-dimensional space $\mathbb U_N \subset \mathbb U^{N_h}$ spanned by the snapshots $u_h(\mu)$ evaluated for the corresponding values of $\mu \in \mathcal P_h$, a finite subset of $\mathcal P$. In fact, once $\mathbb U_N$ is built during the offline phase, a Galerkin projection can be performed for a new value $\mu \in \mathcal P$ to find a surrogate solution $u_N = u_N(\mu) \in \mathbb U_N$
such that
\begin{equation}
    \label{eq:weak_pde_ROM}
        g(u_N,q; \mu) = 
            \left \langle G(u_N; \mu), q \right \rangle = 0 \quad \text{for all } q \in \mathbb U_{N}.
\end{equation}
Among the possible strategies to build the reduced space $\mathbb U_N$, we mention the POD \cite{hesthavenCertifiedReducedBasis2015} and the greedy algorithm \cite{QuarteroniReducedBasisMethods2016,buffa2012priori}. While the POD approach 
compresses a predefined solution dataset, the greedy algorithm adaptively builds $\mathbb U_N$ enriching the reduced space with a properly chosen snapshot at each iteration.
In the following, we denote by $\mathsf B \in \mathbb R^{{N_h}\times N}$
the basis function matrix spanning $\mathbb U_N$, which collects the $N$ basis functions column-wise, and encodes the change of variable from the HF to the RB coordinates. 
Once the basis matrix $\mathsf B$ is provided, the reduced nonlinear system can be solved through the Newton method, whose $k$-th iteration reads: for a given $\mu \in \mathcal P$ and initial guess $\mathsf u_N^{(0)} \in \mathbb R^N$, find $\delta \mathsf u_N \in \mathbb R^N$ such that 
\begin{equation}
\label{eq:Newton_r}
\mathsf {Jac}^G_N(\mathsf u_N^{(k)}; \mu) \delta \mathsf u_N = - \mathsf G_N(\mathsf u_N^{(k)}; \mu), \quad \text{and} \quad \mathsf u_N^{(k+1)} = \mathsf u_N^{(k)} + \delta \mathsf u_N
\end{equation}
where the reduced residual and Jacobian are respectively given by
$$\mathsf G_N(\mathsf u_N^{(k)}; \mu) = \mathsf B^T\mathsf G(\mathsf B\mathsf u_N^{(k)}; \mu), \quad \text{and} \quad \mathsf {Jac}^G_N(\mathsf u_N^{(k)}; \mu) = \mathsf B^T\mathsf {Jac}^G(\mathsf B \mathsf u_N^{(k)}; \mu)\mathsf B.$$
Motivated by the need for an efficient investigation of the parametric space, in this work, we focus on the greedy algorithm, which we now briefly introduce. 

Let $e(\mu) = u_h(\mu) - u_N(\mu)$ be the error between the HF and the reduced solutions for $\mu \in \mathcal P$, we assume to be provided with an error estimator $\Delta_N(\mu)$ such that 
\begin{equation}
\label{estimation}
\lvert \lvert {e}(\mu) \rvert \rvert _{\mathbb U^{N_h}} \leq \Delta_N(\mu).
\end{equation}
Given a tolerance $\varepsilon > 0$ and an initial reduced space $\mathbb U_N = \{u_h( \mu_0)\}$ for \B{a} $\mu_0 \in \mathcal P_h$, the $n$--th step of the algorithm consists in solving the following optimization problem 
$$
 \mu_{n} = \argmax_{\mu \in \mathcal P_h} \Delta_N( \mu),
$$
finding the parameter that maximizes the estimator, and enriching the space with the corresponding snapshot, i.e.,
$\mathbb U_N = \text{span}\{u_h(\mu_0), \dots, u_h(\mu_n)\}$.  At each step, a Gram-Schmidt orthonormalization, indicated by $\mathsf{GS}(\mathsf B, u_h(\mu_n))$, is applied to guarantee the orthogonality of the basis functions $\xi_n$. The approach, from now on referred to as vanilla-greedy and summarized in Algorithm \ref{alg:greedy}, is repeated until a parameter $\mu_n$ verifies
$\Delta_N(\mu_n) \leq \varepsilon$ or the maximum number of basis functions $N_{\text{max}}$ is reached.

\begin{algorithm}
\caption{The vanilla-greedy algorithm}\label{alg:greedy}
\begin{algorithmic}[1]
\Statex{\textbf{Input:} maximum basis number $N_{\text{max}}$, tolerance $\varepsilon$,}
\LineComment{parameter space $\mathcal P_h$, initial parameter $\mu_0$}
\Statex{\textbf{Output:} basis matrix $\mathsf B$} 
\Statex
\State{$\mathsf B = [u_h(\mu_0)] \in \mathbb R^{N_h \times 1}$} \Comment{Initialize the basis function matrix}
\While{$n \leq N_{\text{max}}$ and $\Delta_N(\mu_n) > \varepsilon$}
\State{$\mu_n = \argmax_{\mu \in \mathcal P_h} \Delta_N(\mu)$} \Comment{Find the new parameter $\mu$}
\State{$\xi_n = \mathsf{GS}(\mathsf B, u_h(\mu_n))$} \Comment{Gram-Schmidt orthonormalization}
\State{$\mathsf B \leftarrow [\mathsf B , \xi_n]$} \Comment{Update the basis function matrix}
\EndWhile
\end{algorithmic}
\end{algorithm}

\subsection{A-posteriori error estimator for nonlinear problems}
\label{sec:estimator}
The role of the bound $\Delta_N(\mu)$ is pivotal for the construction of the reduced space during the greedy procedure.

To start our analysis, we consider the \emph{a-posteriori} error estimator for nonlinear PDEs based on the Brezzi-Rappaz-Raviart theory \cite{BRR1,QuarteroniReducedBasisMethods2016}.
By defining the \emph{reduced inf-sup constant} $\beta^h_N (\mu)$ in a neighborhood of the HF representation of the reduced vector $\mathsf u_N$ as 
\begin{equation}
\label{eq:inf_sup_N}
\beta^h_N (\mu) = \adjustlimits \inf_{\mathsf{u} \neq 0} \sup_{\mathsf{q} \neq 0} \frac{\mathsf{q}^T \mathsf{Jac}^G(\mathsf B{\mathsf u_N}; \mu)\ \mathsf{u}}{\norm{\mathsf{u}}_{\mbb U^{N_h}}\norm{\mathsf{q}}_{\mbb U^{{N}_h}}}\qquad \forall \, \mu \in \mathcal P, 
\end{equation}
for $\mathsf u, \mathsf q \in \mathbb R^{N_h}$, the problem is \textit{inf-sup stable} if there exists a positive constant $\hat{\beta}^h_N$ such that 
\begin{equation*}
\beta^h_N (\mu) \geq \hat{\beta}^h_N > 0.
\end{equation*}
The discussion on the role of such a constant is postponed to Section \ref{sec:sota}.

Assuming $ \mathsf{Jac}^G(\mathsf B{\mathsf u_N}; \mu)$ locally Lipschitz in $\mathsf u_N$ with constant $K^h_N (\mu)$, and defining 
\begin{equation}
    \label{eq:tau_general} 
    \tau_N(\mu) = \frac{2K^h_N (\mu) \norm{\mathsf G( \mathsf B \mathsf u_N; \mu)}_{\mathbb U^{{N}_h}}}{\beta^h_N(\mu)^2},
\end{equation}
the \emph{nonlinear error estimator} for nonlinear PDEs is given by
\begin{equation}
\label{eq:nonlinear_est}
\Delta^{\text{nl}}_N(\mu) = \frac{\beta^h_N(\mu)}{K^h_N (\mu)}\left (1 - \sqrt{1 - \tau_N(\mu)} \right ).
\end{equation}
It is clear that the nonlinear estimator is reliable and exploitable only when the condition $\tau_N(\mu)\leq1$ is met for all $\mu \in \mathcal P_h$.
A possible workaround is to employ, even for nonlinear PDEs, the classical \emph{linear bound} 
\begin{equation}
\label{eq:linear_est}
\Delta_N^{\text{lin}}(\mu) = \frac{\norm{\mathsf G( \mathsf B \mathsf u_N; \mu)}_{\mathbb U^{{N}_h}}}{\beta^h_N(\mu)},
\end{equation}
until the nonlinear one can be exploited  \cite{ManzoniEfficientComputationalFramework2014}. We remark that for some problems, e.g., the bifurcating ones, such a condition could be difficult to obtain, and the linear estimator is exploited for all greedy iterations, see 
\ref{app:estimator}.

\section{State of the art for bifurcating phenomena}\label{sec:sota}
In this section, we discuss bifurcating systems and numerical approaches to tackle them. When dealing with these complex scenarios, there are two main problematic features to take into account: (i) the non-uniqueness regime with multiple admissible states, and (ii) the non-differentiability of the parameter-to-solution map at the critical value.
We start by describing the strategies employed in the literature to discover and approximate several coexisting solutions of a bifurcating PDE in the HF setting. Then, we propose a survey of different reduced order approaches for bifurcating problems, to reduce the computational complexity of their analysis.

\subsection{Bifurcation and stability analysis}
\begin{figure}[bht]
    \centering
    \includegraphics{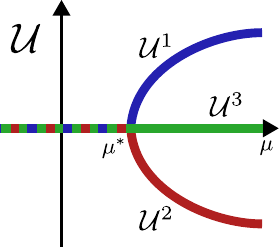}
    \caption{Schematic representation of the solution ensemble in the case of a pitchfork bifurcation.}
    \label{fig:ensamble}
\end{figure}

From the theoretical standpoint, the existence and uniqueness of the solution to nonlinear parametric PDEs are guaranteed by the implicit function theorem \cite{ambrosetti1995primer,ciarlet2013linear}.
When the regularity assumption is satisfied, the solution evolves uniquely and continuously w.r.t.\ the parametric dependence. In contrast, when a small variation of the parameter under consideration causes a drastic change in the system's response, the problem admits coexisting solutions with a change in their stability properties, undergoing a \textit{bifurcating phenomenon}. The value of the parameter corresponding to the loss of uniqueness is referred to as \emph{bifurcation point}, and we denote it by $\mu^*$.\
The bifurcation is due to the non-invertibility of the Fr\'echet derivative of $G$ at the bifurcation point, which represents the failure of the implicit function theorem.
In other words, when a bifurcation occurs, $G$ is no longer an injective map around $u(\mu^*)$, and admits multiple configurations. At the discrete level, the local invertibility of $G$ translates into the positivity of the {discrete Babu{\v s}ka inf-sup condition} \eqref{eq:inf_sup_disc}, i.e., $\beta^h (\mu) \geq \hat{\beta}^h > 0$. This condition does not hold for $\mu^*$. For further details, the reader may refer to \cite{QuarteroniReducedBasisMethods2016,pichi_phd}.  
In this scenario, we define the \emph{solution branches} as the multiple physical behaviors admitted by the system while varying $\mu \in \mathcal{P}$.
Assuming the coexistence of up to $K$ branches, meaning $K$ different possible configurations, we denote with $\mathcal U^i$ the $i$-th branch for $\mu \in \mathcal P$, which collects the set of solutions $u^i_h(\mu)$ as the continuation of the solution sharing same qualitative properties. 
The union of all the solutions branches, i.e., the \emph{solution ensemble}, schematically depicted in Figure \ref{fig:ensamble}, is defined as follows:
\begin{equation}
\label{ensemble}
\mathcal U = \{ u^i_h(\mu) \in \Cal U^i \; | \; \mu \in \mathcal P \}_{i = 1}^K .
\end{equation}

\subsection{Deflation method}
\label{sec:def}
An important question in bifurcation analysis is how to identify and fully capture all the many possible coexisting branches that characterize the bifurcating phenomenon.
Towards this goal, we now discuss the basic concept of the \emph{deflation} algorithm for parametric PDEs \cite{Farrell2015A2026}. This strategy aims at modifying the PDE residual to ``penalize'' the convergence towards an already known solution and ``drive'' the system to discover a new one. From the mathematical standpoint, given a solution $u^{j}_h(\underline{\mu}) \in \Cal U^j$, for fixed ${j}$ and $\underline{\mu} \in \mathcal P$, the algorithm deflates the nonlinear system such that the Newton method converges, in the bifurcating regime, to a solution $u^{\ell}_h(\underline{\mu})$ with $\ell \neq j$, and diverges elsewhere. 
The algorithm can be iterated over all branches, reaching a full discovery of the coexisting solutions for $\underline{\mu}$, and finally, iterating over $\mathcal P$,
reconstructs the solution ensamble $\mathcal U$.
The deflation method builds a new residual vector $F$ from the residual $G$ defined in \eqref{eq:strong_pde}, such that: (i) a solution for $G$ is also a solution for $F$ and (ii) given a set of known solutions for a fixed $\underline \mu$, the solver will not converge to any of them.

Let us denote by $u^1_h = u^1_h(\underline{\mu})$ the solution of Equation \eqref{eq:strong_pde} for a fixed $\underline{\mu} \in \mathcal P$ with initial guess $u^{(0)}$. 
We can now build the modified residual as a function of the (renamed) unknown $y_h \in \mathbb{U}^{N_h}$
\begin{equation}
\label{eq:F_mod}
    F(y_h; \underline{\mu}) = M(y_h, u^1_h)G(y_h; \underline{\mu}),
\end{equation}
where, given $\mathbb I$ the identity over $\mathbb U$, $r \geq 1$, and $\sigma > 0$, the \textit{deflation operator} $M(y_h, u^1_h)$ is defined as
\begin{equation}
    \label{eq:def_M}
    M(y_h, u^1_h) = \left ( \frac{1}{\norm{y_h - u^1_h}^r} + \sigma \right ) \mathbb I.
\end{equation}

From the definition of the deflated residual in Equation \eqref{eq:F_mod}, it is clear that when the Newton method tries to converge to the already known solution $u^1_h$, the denominator tends to zero and guides the iteration to converge to a newly discovered solution $y_h = u^2_h  = u^2_h(\underline{\mu})$, or to diverge. The deflation process can be performed iteratively until the solver diverges for the system
\begin{equation}
\label{eq:F_mod_final}
    F(y_h; \underline{\mu}) = M(y_h, u^1_h)M(y_h, u^2_h)\cdots M(y_h, u^n_h)G(y_h; \underline{\mu}),
\end{equation}
after the discovery of $n$ distinct coexisting solutions.
We remark that, in general, the deflation method does not guarantee to reach all possible solutions for a given parameter, i.e., $n \leq K$, but insights into their convergence can be obtained through suitable assumptions \cite{RheinboldtAdaptiveContinuationProcess1978,RallNoteConvergenceNewtons1974,FarrellComputationDisconnectedBifurcation2016}.
The $k$-th step of the Newton method applied to the deflated system $F(y_h,{\underline{\mu}})$ reads as: find $\delta \mathsf y_h \in \mathbb R^{N_h}$ such that
\begin{equation}
\label{eq:Newton_method_k_F}
\mathsf {Jac}^F(\mathsf y_h^{(k)}; \underline{\mu}) \delta \mathsf y_h = - \mathsf F(\mathsf y_h^{(k)}; \underline{\mu}).    
\end{equation}
One of the main features of the deflation algorithm is its scalability, meaning that, thanks to the Shermann-Morrison formula (for further details see \cite{Farrell2015A2026}), one can find a solution to Equation \eqref{eq:Newton_method_k_F} exploiting the original problem \eqref{eq:strong_pde} as
\begin{equation}
    \label{eq:tau_system}
    \delta \mathsf y_h = \theta\delta \mathsf u_h, 
    \qquad
    \text{with}
    \qquad
    \theta = 1 - \frac{\mathsf M^{-1} {(\mathsf {Jac}^G)}^{-1} \mathsf G \mathsf E^T  }{1 + \mathsf  E^T \mathsf  M^{-1} {(\mathsf {Jac}^G})^{-1} \mathsf G},
\end{equation}
where $\mathsf M$ and $\mathsf E$ are the algebraic representations of the operators $M = \prod_{i=1}^n M(y_h, u^i_h)$ and its derivative $M'$, respectively.
The deflation procedure is described in Algorithm \ref{alg:def-newton}.
\begin{algorithm}
\caption{The $ \mathsf {DeflatedNewton}$ algorithm}\label{alg:def-newton}
\begin{algorithmic}[1]
\Statex{\textbf{Input:} residual vector $\mathsf G$, tolerance $\mathsf {tol}$, initial guess $y^{(0)}_h$, parameter $\mu$}
\Statex{\textbf{Output:} solution of the deflated system $\mathsf y$} 
\Statex
\While{$\norm{\mathsf G(\mathsf y_h^{(k)}; \mu)}_{\mathbb U^{N_h}} \geq \mathsf{tol}$}
\State{$\mathsf {Jac}^G(\mathsf y^{(k)}_h; \mu) \delta \mathsf u_h = - \mathsf G(\mathsf y^{(k)}_h; \mu)$} \Comment{Solve original linearized system}
\State{Compute $\theta$} \Comment{Deflation step}
\State{$\mathsf y^{(k+1)}_h = \mathsf y^{(k)}_h + \theta \delta \mathsf u_h$} \Comment{Update solution}
\State{$k = k + 1$}
\EndWhile
\end{algorithmic}
\end{algorithm}

\begin{remark}[Continuation methods]
When approximating branches of solutions belonging to $\mathcal U$, an important class of algorithms conceived to follow the bifurcating phenomenon is represented by the \emph{continuation} methods. The main idea of the continuation is to provide a close enough initial guess for the nonlinear solver, such that it requires a cheaper computational cost and, most importantly, the new parametric solution shares the same qualitative properties as the guess. There exist many different approaches, from predictor-corrector schemes to arclength strategies \cite{KellerLecturesNumericalMethods1988,UeckerNumericalContinuationBifurcation2021,AllgowerNumericalContinuationMethods1990,LegerImprovedMoorePenroseContinuation2023a}, but for the sake of simplicity here we will only consider a simple continuation method, choosing the last computed solution corresponding to the previous value of the parameter as the initial guess.
\end{remark}

\subsection{ROMs for bifurcation systems}
\label{sec:rombif}
As we have seen in the previous sections, the approximation of the bifurcation diagram usually requires an unbearable computational cost. Indeed, one needs to solve high-dimensional linear systems for every iteration of the Newton method, for each parameter in the discretized parametric space to reconstruct the branch, and for every branch existing in the bifurcation diagram.
For these reasons, in the ROM literature, the POD has been extensively used to compress and extract the main information from the snapshots matrix \cite{pichi_phd,PittonComputationalReductionStrategies2017,HerreroRBReducedBasis2013}, with the final goal of accelerating such an analysis. Different strategies are possible: considering global (i.e., the snapshots are collected over $\mathcal U$) or branch-wise reductions (i.e., a specific reduced model is built for the branch $\mathcal U^i$, with $i=1, \dots, K$) \cite{PichiReducedBasisApproaches2019,PichiReducedOrderModeling2020,PichiDrivingBifurcatingParametrized2022a} or local ROMs \cite{HessLocalizedReducedorderModeling2019,CortesLocalROMRayleigh2024}, but all of them require the \emph{a-priori} information of the bifurcating phenomenon. Indeed, while POD approaches do not need error estimators and are only related to data compression based on energy criteria, no certification of the error is guaranteed, and the reduced solver may converge to a different branch depending on how the reduced space has been built.

In contrast, as we discussed in Section \ref{sec:ROMnonlinear}, the certified greedy approach is not straightforward for nonlinear problems and is challenged by the local invertibility of the system. Indeed, the error certification is strictly related to: 
\begin{enumerate}
\item[(i)] the value of $\beta^h_N(\mu)$ and the well-posedness of the system,
\item[(ii)] a \emph{unique} value of the (linear or nonlinear) estimator $\Delta_N(\mu)$ for a given $\mu \in 	\mathcal P$.
\end{enumerate}
In the case of a bifurcating system, the inf-sup condition is not verified at the bifurcation point $\mu^*$, and both linear and nonlinear estimators cannot be used.
Moreover, when multiple solutions arise, the estimator should provide a set of values $\{\Delta_N^i(\mu)\}_i$ related to the solution on the different branches, and not a unique information. 
As we will see in Section \ref{sec:results}, the numerical results confirm that applying the vanilla-greedy algorithm to bifurcating PDEs can lead to only a partial representation of the system behavior since the estimator lacks reliability and the action of the algorithm is spoiled.
In what follows, we discuss novel greedy algorithms that are tailored for bifurcating systems and conceived to overcome these challenges.

\section{The adaptive-greedy algorithm}
\label{sec:adaptivegreedy}
We now propose an enhanced version of the vanilla-greedy, summarized in Algorithm \ref{alg:ref-greedy}: the adaptive-greedy. 
The main purpose of this strategy is to build an adaptive training set related to the detection of the bifurcation points of a system based on successive approximations of the reduced inf-sup constant \eqref{eq:inf_sup_N}. In fact, since the assumption in Equation \eqref{eq:inf_sup_disc} is not verified at the bifurcation point, numerically, this feature translates into observing small values of $\beta^h(\mu)$, for $\mu$ close enough to $\mu^*$. Thus, if the parametric space is sampled properly, we can assume that $\mu^*$ is well approximated by the parameter that minimizes the discrete inf-sup constant.  At the reduced level, for large enough $N$, the reduced inf-sup \eqref{eq:inf_sup_N} is a good approximation of the high-fidelity inf-sup \eqref{eq:inf_sup_disc}. This is the rationale for approximating the bifurcation point as the one for which $\beta^h_N(\mu)$ is minimized. 

\begin{algorithm}
\caption{The adaptive-greedy}\label{alg:ref-greedy}
\begin{algorithmic}[1]
\Statex{\textbf{Input:} maximum iteration $N_{\text{max}}$, tolerance $\varepsilon$}
\LineComment{an \emph{initial} parameter space $\mathcal P_h$, initial parameter $\mu_0$.}
\Statex{\textbf{Output:} basis matrix $\mathsf{B}$, bifurcation point $\mu^*$} 
\Statex
\State{$\mathsf B = [u_h(\mu_0)] \in \mathbb R^{N_h \times 1}$} \Comment{Initialize the basis function matrix}
\While{$n \leq N_{\text{max}}$ and $\Delta_N(\mu_n) > \varepsilon$}
\State{$\mu_n = \argmax_{\mu \in \mathcal P_h} \Delta_N(\mu)$} \Comment{Find the new parameter $\mu$}
\State{$\xi_n = \mathsf{GS}(\mathsf B, u_h(\mu_n))$} \Comment{Gram-Schmidt orthonormalization}
\State{$\mathsf B \leftarrow [\mathsf B , \xi_n]$} \Comment{Update the basis function matrix}
\State{$\mu_{\text{bif}} = \argmin_{\mu \in \mathcal P_h} \beta_N^h(\mu)$} \Comment{Approximating the bifurcation point}
\If{$n > 1$}
\State{$\mathcal P_h = \mathsf {Refinement}(\mathcal P_h, \mu_{\text{bif}}, \mu_{\text{prev}}, n_{\text{ref}}, \mathsf {tol})$} \Comment{Update training set and $\mu_{\text{bif}}$}
\EndIf
\State{$\mu_{\text{prev}} = \mu_{\text{bif}}$}
\State{$ n =  n + 1$} 
\EndWhile
\State{$\mu^* = \mu_\text{guess}$}
\end{algorithmic}
\end{algorithm}

The procedure adds points to the training set $\mathcal P_h$ in a neighborhood of the current approximation of the bifurcation point and, lastly, detects $\mu^*$ itself.
We now describe the procedure for the detection of a single bifurcation point combined with the vanilla-greedy strategy. We compute the parameter that minimizes the reduced inf-sup constant of $\mathcal P_h$ for $n=1$ after the first basis enrichment, i.e., when $\mathsf B \in \mathbb R^{N_h \times 2}$, and we denote it as $\mu_{\text{bif}}$. 

In the next iteration, we enrich the basis with the standard greedy policy, we store the previous approximation $\mu_{\text{bif}}$ as $\mu_{\text{prev}}$, and we compute a new $\mu_{\text{bif}}$ minimizing the reduced inf-sup constant. We now can apply the \textsf{Refinement} function, described in Algorithm \ref{alg:refine}, where we compare $\mu_\text{prev}$ and $\mu_\text{bif}$: if their distance is greater than a prescribed tolerance, a refinement of $n_{\text{ref}}$ points is performed around $\mu_{\text{bif}}$. The training set is finally enriched with these new parameters.

In our numerical experiments, we refine the interval $(\mu_a, \mu_b)$, where the extremes represent the left and right neighbors of $\mu_{\text{bif}}$ in $\mathcal P_h$, which we suppose to be ordered without loss of generality. In case $\mu_{\text{bif}}$ is the first element or the last one of $\mathcal P_h$, we refine the intervals $(\mu_{\text{bif}}, \mu_b)$ and $(\mu_a, \mu_{\text{bif}})$, respectively. For the sake of brevity, we do not report these extreme cases in Algorithm \ref{alg:refine}. 

The procedure is repeated iteratively until $\mu_{\text{bif}}$ and $\mu_{\text{prev}}$ meet the similarity criterion. At the end of the procedure, $\mu_{\text{bif}}$ represents a reliable approximation of the critical point. 
Finally, we recall that this procedure is completely independent of the deflation.

\begin{algorithm}
\caption{The $ \mathsf {Refinement}$ function}\label{alg:refine}
\begin{algorithmic}[1]
\Statex{\textbf{Input:} parameter space $\mathcal P_h$, bifurcation point guesses $\mu_{\text{bif}}$ and $\mu_{\text{prev}}$,}
\LineComment{number of points to add $n_{\text{ref}}$, tolerance $\mathsf {tol}$}
\Statex{\textbf{Output:} new training set $\mathcal P_h$}
\Statex
\If {$|\mu_{\text{bif}} - \mu_{\text{prev}}| > \mathsf{tol}$}
 \State{$\mu_a$ = the left neighbor of $\mu_{\text{bif}}$ in $\mathcal P_h$}
 \State{$\mu_b$ = the right neighbor of $\mu_{\text{bif}}$ in $\mathcal P_h$}
 
\State{Sample $\{\overline \mu_{i}\}_{i=1}^{n_{\text{ref}}} \in (\mu_a, \mu_b)$}
\State{$\mathcal P_h = \mathcal P_h \cup \{\overline \mu_{i}\}_{i=1}^{n_{\text{ref}}}$}
\EndIf 
\end{algorithmic}
\end{algorithm}

The adaptive-greedy has the following advantages w.r.t.\ the vanilla-greedy:
\begin{itemize}
\item[$\circ$] the adaptive strategy allows one to start with a coarse $\mathcal P_h$, concentrate the sample around the bifurcating point that is eventually detected, and save computational effort avoiding the exploration of the uniqueness regime, which by definition requires less information to be represented;
\item[$\circ$] the algorithm can be seen as a pre-processing step for the deflated-based strategies, which can be applied once the bifurcation point has been detected. In other words, the deflation can be performed only on the subset of $\mathcal P_h$ featuring multiple solutions, reducing the time required to construct the basis functions.
\end{itemize}

\section{The deflated-greedy algorithm}
\label{sec:defgreedy}
\begin{figure}[htb]
	\includegraphics[height=0.23\textheight]{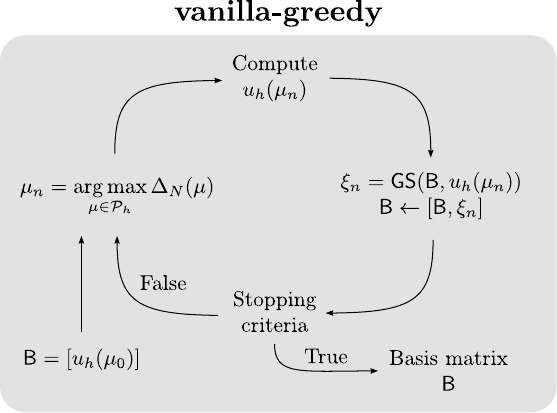}\hfill
	\includegraphics[height=0.23\textheight]{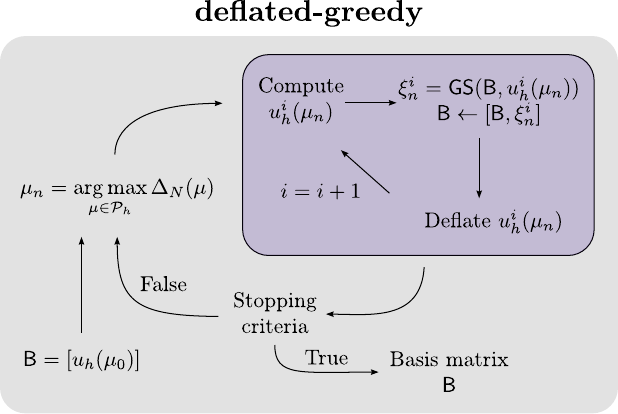}
	\caption{Comparison of vanilla-greedy methodology and the proposed deflated-greedy approach.}
    \label{fig:greedyschemes}
\end{figure}

In this section, we present a novel deflation-based Greedy algorithm, which aims to provide error certification for all branches of a bifurcating system, employing the deflation and continuation methods introduced in Section \ref{sec:sota} without any \emph{a-priori} information on the bifurcating phenomenon itself. We illustrate the proposed methodology and its main differences w.r.t.\ the vanilla approach in Figure \ref{fig:greedyschemes}.

\subsection{Reduced deflation method}
\label{sec:reduced_def}

To obtain a greedy strategy that is capable of certifying the whole set of (non-unique) solutions $\mathcal U$ for general nonlinear systems we need to transpose the concept of deflation and continuation at the reduced level. While the latter is straightforward, simply entailing the same mechanism but in a lower dimensional space, we discuss the \emph{reduced deflation algorithm} exploited in the proposed methodology below. In analogy to the HF approach, we aim at approximating the \emph{reduced solutions set} defined as
\begin{equation}
\label{ensembleROM}
\mathcal U_N = \{ u^i_N(\bmu) \in \Cal U_N^i \in \mathbb{R}^N \; | \; \bmu \in \mathcal P \}_{i = 1}^K .
\end{equation}
Fixed a parameter $\underline{\mu} \in \mathcal P$, we apply the reduced Newton algorithm in Equation \eqref{eq:Newton_r} with a suitable initial guess, obtaining a first reduced approximated solution $u^1_N = u^1_N(\underline{\bmu})$. Then, we deflate the original system w.r.t.\ $u^1_N$ and look for other admissible solutions belonging to different reduced branches. To do so, at the $n$-th step, we aim at solving the reduced deflated system given by the residual 
\begin{equation}
\label{eq:F_mod_r}
    F(y_N; \underline{\bmu}) =  M(y_N, u^1_N)M(y_N, u^2_N)\cdots M(y_N, u^n_N)G(y_N; \underline{\bmu}).
\end{equation}
Instead of applying the Newton method to the deflated residual, we solve Equation \eqref{eq:Newton_r} around $\mathsf y_N$, and solve for the modified update 
\begin{equation}
\label{eq:tau_N}
\delta \mathsf y_N = \theta_N \delta \mathsf u_N,
\qquad
\text{with}
\qquad
    \theta = 1 - \frac{\mathsf M^{-1}_N {(\mathsf {Jac}^G_N)}^{-1} \mathsf G_N \mathsf E^T_N  }{1 + \mathsf  E^T_N \mathsf  M^{-1}_N {(\mathsf {Jac}^G_N})^{-1} \mathsf G_N},
\end{equation}
where $\mathsf M_N$ and $\mathsf E_N$ are the reduced matrix representations of the operators $M_N = \prod_{i=1}^n M(y_N, u^i_N)$ and $M'_N$, respectively. 
\B{Thus, we have developed a reduced counterpart of the high-fidelity deflation framework presented in Section \ref{sec:def}. Such an approach will be used to find coexisting ``reduced" branches in the low-fidelity approximations needed for both the error certification and the actual reconstruction, performed during the offline and online phases, respectively. In particular, given a novel parameter in the continuation setting, previous solutions are exploited: (i) to follow the evolution of the branches and (ii) to deflate already known states to find novel ``reduced" branches.}

\subsection{Certification of multiple branches} 
\label{sec:alg1}

The deflated-greedy methodology is a tailored greedy strategy to obtain certified ROMs that combines the reduced deflated method and the iterative greedy algorithm to build the reduced space. It presents two main differences w.r.t.\ the vanilla-greedy algorithm described in Section \ref{sec:greedy}:
\begin{itemize}
\item[$\circ$] multiple snapshots corresponding to a same value of the parameter might be picked via HF deflation;
\item[$\circ$] multiple evaluations of the (deflated) error estimator drive the search towards the solution branch and parametric regions less represented by the current basis functions.
\end{itemize}
In what follows, we describe the deflated-greedy procedure summarized in Algorithm \ref{alg:defgreedy}.

\begin{algorithm}
	\caption{The deflated-greedy algorithm}\label{alg:defgreedy}
	\begin{algorithmic}[1]
	\Statex{\textbf{Input:} maximum number of basis function $N_{\text{max}}$, tolerance $\varepsilon$,}
	\LineComment{training set $\mathcal P_h$, initial parameter $\mu_0$, HF guesses set $\Xi_{HF}$,}
	\LineComment{RB guesses set $\Xi_{RB}$, tolerance $\text{tol}$}
	\Statex{\textbf{Output:} final basis function matrix $\mathsf B$}
	\Statex
	\State{$\mathsf B = [u_h(\mu_0)] \in \mathbb R^{N_h \times 1}$}, $n = 1$ \Comment{Initialize the basis function matrix}
	\State{$\Delta = [\Delta(\mu_0)]$}\Comment{Initialize the estimator set}
	\While{$n \leq N_{\text{max}}$ and $\max \{\Delta\} > \varepsilon$}
	\State{$\Upsilon_{HF} = []$}
	\State{$\Delta = \mathsf{DeflatedEstimator}(\mathcal P_h, \Xi_{RB}, \mathsf B, \text{tol})$}\Comment{Compute the deflated estimator}
	\State{$\mu_n = \argmax_{\mu \in \mathcal P_h}{\Delta}$}
	\Comment{Find the new $\mu$}
	\State{${u_N} = \argmax_{u_N \in \mathcal U_N^i} \{ \Delta^i_N(\mu_n)\}_{i=1}^{b_{\mu_n}}$}\Comment{Find the worst approximated branch}
	\State{$u_h(\mu_n) = \mathsf{Newton}(\mathsf G, \text{tol}, \mathsf B \mathsf u_N)$} \Comment{New snapshot}
	\State{$\xi_n = \mathsf{GS}(\mathsf B, u_h(\mu_n))$} \Comment{Gram-Schmidt orthonormalization}
	\State{$\mathsf B \leftarrow [\mathsf B , \xi_n]$} \Comment{Update the basis function matrix}
	\State{$\Upsilon_{HF} \leftarrow [\Upsilon_{HF} , u_h(\mu_n)]$} \Comment{Update the roots}
	\State{$\mathsf B, \Xi_{HF}, n =  \mathsf{DeflatedSnapshots}(\Upsilon_{HF}, \Xi_{HF},\mu_n, \text{tol}, \mathsf B)$}\Comment{Enrich the basis with deflated snapshots}
	\EndWhile
	\end{algorithmic}
\end{algorithm}

We start by setting as a first basis the initial solution $u_h(\mu_0)$, for $\mu_0$ in $\mathcal P_h \subset \mathcal P$. Then, we compute the error estimator using the \textsf{DeflatedEstimator} function (see Algorithm \ref{alg:defestimator}). 

For the first parameter of the ordered set $\mathcal P_h$, we solve the reduced problem in Equation \eqref{eq:weak_pde_ROM} with a trivial initial guess. This way, a first reduced solution is found, populating the set of reduced roots denoted with $\Upsilon_{RB}$. Then, we proceed with the \textsf{DeflatedNewton} strategy, i.e., we apply Algorithm \ref{alg:def-newton} to the reduced residual $\mathsf G_N$.

Since we are interested in using different initial conditions for the reduced Newton solver, we denote the set of reduced initial conditions by $\Xi_{RB}$. This set is intended as a way to propagate the information by continuing the discovered branches.

We populate the set $\Upsilon_{RB}$ of reduced roots until the algorithm reaches the termination condition, and we compute the error estimator $\Delta_N(\mu)$ corresponding to the solutions we found. In fact, assuming that for $\mu \in \mathcal P_h$ we found $b_{\mu}$ solutions, we compute the estimator for each of them, leading to the \emph{estimator set} $\Delta = \{\Delta^i_N(\mu) \; | \; \mu \in \mathcal P_h\}_{i=1}^{b_{\mu}}$. 
It is clear that the number of discovered reduced solutions $b_{\mu}$ changes w.r.t.\ $\mu$ due to their existence and/or convergence issues, and might be lower than the number of branches $K$.
Before moving to the next parameter, we define the set of reduced guesses $\Xi_{RB}$ as the set of computed reduced roots. Namely, we exploit the continuation approach with the reduced solutions of the previous step to drive the exploration in the direction of the obtained branches.

At this point, we are interested in the parameter $\mu_n$, corresponding to the reduced solution $u_N^i(\mu_n)$, that maximizes the estimator set $\Delta$, seeking the worst approximated branch.

Then we compute a new solution $u_h(\mu_n)$ with the reduced solution projected onto the HF basis as an initial guess for the Newton solver, to capture the information of the less represented branch.
After a $\mathsf {GS}$ orthonormalization, the solution $u_h(\mu_n)$ is used to enrich the reduced space and the set of HF roots, defined as $\Upsilon_{HF}$.

We now aim to discover new branches by applying the standard deflation-based approach to the HF system by means of the \textsf{DeflatedSnapshots} function (see Algorithm \ref{alg:defsnap}). 
The deflation algorithm is applied for several initial guesses stored in the set $\Xi_{HF}$, including all the basis functions computed so far, and thus encoding all the information from the previous iterations. 
\begin{algorithm}
\caption{The $\mathsf{DeflatedEstimator}$ function}\label{alg:defestimator}
\begin{algorithmic}[1]
\Statex{\textbf{Input:} training set $\mathcal P_h$, RB guesses set $\Xi_{RB}$,}
\LineComment{basis function matrix $\mathsf B$, tolerance $\text{tol}$}
\Statex{\textbf{Output:} $\Delta = \{\Delta^i_N(\mu)\;| \;\mu \in \mathcal P_h \}_{i=1}^{N_b^{\mu}}$}
\Statex
\State{$\Delta = []$}\Comment{Initialize the estimator set}
\For{$\mu \in \mathcal P_h$}
\State{$\Upsilon_{RB} = []$}\Comment{Initialize the reduced roots}
\State{$u_N(\mu) = \mathsf{Newton}(\mathsf G_N, \text{tol})$}\Comment{Solve with continuation guess}
\State{$\Upsilon_{RB} \leftarrow [\Upsilon_{RB}, u_N(\mu)]$} \Comment{Update reduced roots}
\For{$v_N \in \Xi_{RB}$}
\While{reduced deflation converges}
	\State{$M_N = \prod_{u_r \in \Upsilon_{RB}} M(y_N(\mu), u_r)$}\Comment{Assemble reduced deflation operator}
	\State{$y_N(\mu) = \mathsf{DeflatedNewton}(\mathsf G_N, \mathsf M_N, \text{tol}, \mathsf B, \mathsf v_N)$}\Comment{Deflated reduced solve}
\State{$\Upsilon_{RB} \leftarrow [\Upsilon_{RB}, y_N(\mu)]$} \Comment{Update reduced roots}
\State{$ \Delta \leftarrow [\Delta, \Delta_N(\mu)]$} \Comment{Update estimator set}
\EndWhile
\EndFor
\State{$\Xi_{RB} = \Upsilon_{RB}$} \Comment{Exploit reduced roots as guesses set}
\EndFor
\end{algorithmic}
\end{algorithm}
Any new deflated solution enriches the basis function matrix $\mathsf B$ and the set of HF roots $\Upsilon_{HF}$, defining $n$ as the total number of basis functions.

The deflated-greedy procedure ends when a maximum number of basis functions $N_{max}$ is reached or a convergence criterion over the estimator set $\Delta$ is met.

\begin{remark}[Error estimator]
\label{rem:check}
The same parameter can be picked in two different iterations of the deflated-greedy. Indeed, while a basis can easily represent a branch, it can be useless for a coexisting solution featuring qualitatively different physical phenomena (practical examples will be provided in Section \ref{sec:results}). Namely, the same parameter can maximize the estimator multiple times throughout the deflated-greedy iterations. In this case, the HF solution of line 6 of Algorithm \ref{alg:defgreedy}, even if guided toward the behavior that is less represented by means of the initial guess, can converge to a solution that is already in the set of basis functions.
For these reasons, we need to check each new deflated solution to avoid the addition of redundant information.
If no new information has been added to the basis, we consider the second-last parameter that maximizes $\Delta$, and we proceed until a new basis is discovered. 
\end{remark}
\begin{algorithm}
\caption{The $\mathsf{DeflatedSnapshots}$ function}\label{alg:defsnap}
\begin{algorithmic}[1]
\Statex{\textbf{Input:} HF roots set $\Upsilon_{HF}$, HF guesses $\Xi_{HF}$, a parameter $\mu$,}
\LineComment{tolerance tol, basis function matrix $\mathsf B$}
\Statex{\textbf{Output:} the updated basis matrix $\mathsf B,$ the updated HF guess list $\Xi_{HF}$,}
\LineCommentO{the total number of basis functions $n$}
\Statex
\For{$u \in \Xi_{HF}$}
\While{deflation converges}
	\State{$M = \prod_{u_r \in \Upsilon_{HF}} M(y_h(\mu), u_r)$}\Comment{Assemble deflation operator}
	\State{$y_h(\mu) = \mathsf{DeflatedNewton}(\mathsf G, \mathsf M, \text{tol}, \mathsf B, \mathsf u)$}\Comment{Deflated solve}
\State{$\Upsilon_{HF} \leftarrow [\Upsilon_{HF} , y_h(\mu))]$} \Comment{Update the roots}
	\State{$\xi_n = \mathsf{GS}(\mathsf B, y_h(\mu))$}\Comment{Gram-Schmidt orthonormalization}
	\State{$\mathsf B \leftarrow [\mathsf B , \xi_n]$} \Comment{Update the basis function matrix}
 \State{$n =$ number of basis functions}
\EndWhile
\EndFor

\State{$\Xi_{HF} \leftarrow [\Xi_{HF}, \Upsilon_{HF}]$} \Comment{Update roots}
\end{algorithmic}
\end{algorithm}

\A{
\begin{remark}[Computational times] As anticipated in Section \ref{sec:ROMnonlinear}, in this context, we are not focusing on fully retrieving the efficiency, as needed when dealing with nonlinear problems. The main reason is that, by introducing an approximation of the nonlinear term via hyper-reduction approaches, the complex bifurcating structure could be completely lost and/or create non-converge issues, as already noticed in the literature \cite{PichiReducedOrderModels2023,PichiReducedOrderModeling2020}.
Moreover, concerning the computational cost of augmenting the greedy strategy by deflation, providing a precise and complete analysis is quite challenging since it is very much dependent on the problem's bifurcating structure.
In fact, the vanilla-greedy strategy requires a single high-fidelity computation for each iteration, but at the same time, it limits the number of branches that can be recovered. The iterations of the deflated-greedy strategy are, in general, computationally more costly but allow for the potential discovery of novel branches. Moreover, even exploiting some continuation strategy, the vanilla approach could stall sampling the same branch over and over, with no approximation capability on the other branches. Thus, the deflated variant always provides more information and could also avoid computational bottlenecks, reducing the time required for the low-fidelity approximation.
\end{remark}
}

\B{
\begin{remark}[Bifurcation agnostic] 
We highlight that the developed framework is meant to exploit no prior information on the position, number, and type of bifurcating phenomena. Indeed, while the adaptive procedure, driven by the reduced model discovery itself, takes care of an efficient parametric sampling and bifurcation points identification, the combined exploitation of greedy error estimator and deflation strategy allows the approximation of several coexisting branches both at the high-fidelity and reduced-order levels.
\end{remark}
}

\section{Numerical results}
\label{sec:results}
In this section, we show the advantages of exploiting the adaptive- and deflated-greedy algorithms in bifurcating contexts. 
We test the two novel ROM approaches on a channel flow model undergoing a pitchfork bifurcating phenomenon, i.e., the Coanda effect. The obtained results confirm the capability of the novel algorithms to outperform POD and the vanilla-greedy approach when dealing with non-uniqueness behaviors, providing a complete approximation and efficient detection of the possibly coexisting states. 

The computations in this work have been performed with RBniCS \cite{rbnics} library, which is an implementation in FEniCS \cite{fenics} of several reduced order modeling techniques.

\subsection{Numerical setting: the Coanda effect}
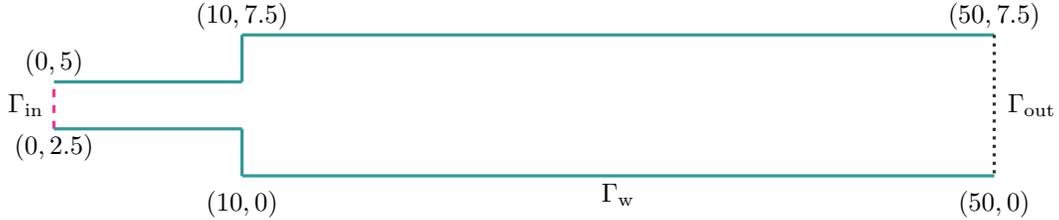
\begin{figure}[h]
\begin{center}
\begin{tikzpicture}[scale=.25]

\filldraw[color=magenta!90, very thick, dashed](0,2.5) -- (0.,5.);
\filldraw[color=teal!80, fill=gray!10, very thick](0,2.5) -- (10,2.5);
\filldraw[color=teal!80, fill=gray!10, very thick](0,5) -- (10,5);
\filldraw[color=teal!80, fill=gray!10, very thick](10,5) -- (10,7.5);
\filldraw[color=teal!80, fill=gray!10, very thick](10,2.5) -- (10,0);
\filldraw[color=teal!80, fill=gray!10, very thick](10,0) -- (50,0);
\filldraw[color=teal!80, fill=gray!10, very thick](10,7.5) -- (50,7.5);
\filldraw[color=black!80, fill=gray!10, very thick, dotted](50,7.5) -- (50,0);

\node at (-1.5,3.75){\color{black}{$\Gamma_{\text{in}}$}};
\node at (52,3.75){\color{black}{$\Gamma_{\text{out}}$}};
\node at (30,-1){\color{black}{$\Gamma_{\text{w}}$}};
\node at (0,1.5){\color{black}{$(0,2.5)$}};
\node at (0,6){\color{black}{$(0,5)$}};
\node at (10, -1.5){\color{black}{$(10,0)$}};
\node at (50, -1.5){\color{black}{$(50,0)$}};
\node at (50, 8.5){\color{black}{$(50,7.5)$}};
\node at (10, 8.5){\color{black}{$(10,7.5)$}};

\end{tikzpicture}
\end{center}
\caption{Sudden-expansion channel domain $\Omega$.}
\label{fig:channel}
\end{figure}

The Coanda effect is a well-known bifurcating phenomenon modelled by the Navier-Stokes equations in a sudden expansion channel, and describes the symmetry-breaking behavior and the stability properties of a fluid flow corresponding to different values for the kinematic viscosity \cite{tritton2012physical}. The two-dimensional computational domain, depicted in Figure \ref{fig:channel}, is denoted with $\Omega$.

\begin{figure}
    \centering
\begin{minipage}{0.45\textwidth}
    \includegraphics[width=\linewidth]{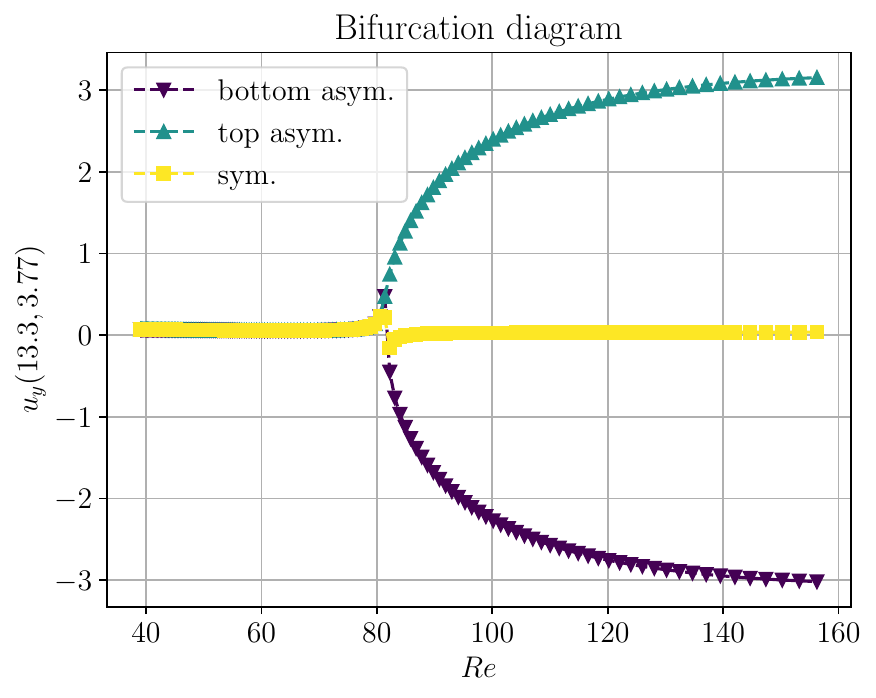}
\end{minipage}\hfill
\begin{minipage}{0.5\textwidth}
    \includegraphics[width=\linewidth]{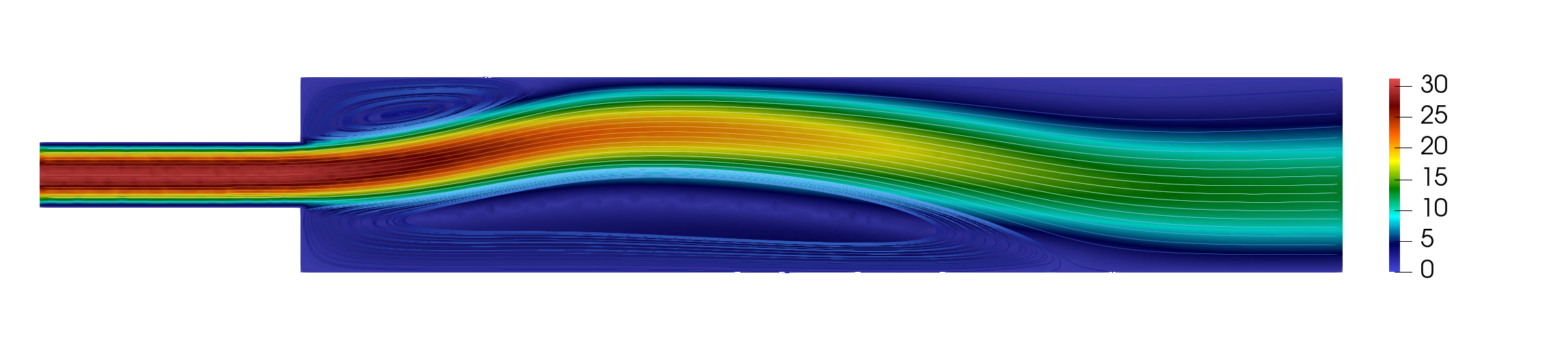}    
    
    \includegraphics[width=\linewidth]{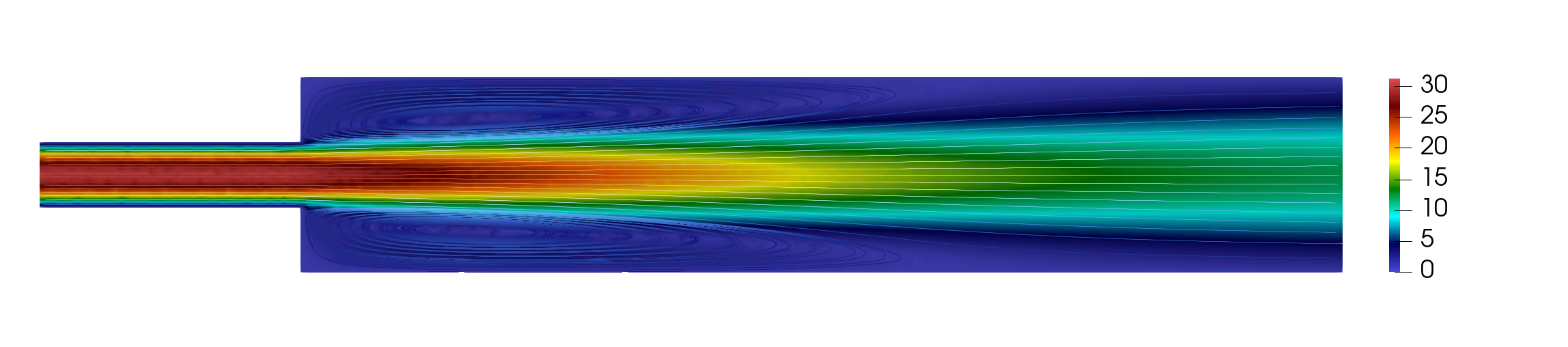}    
    
    \includegraphics[width=\linewidth]{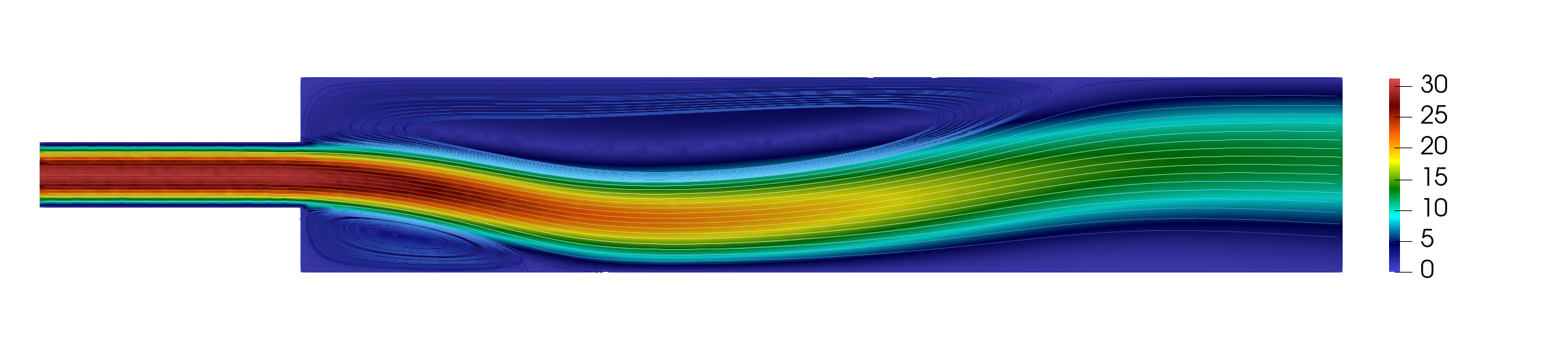}
        \end{minipage}
    \caption{Pitchfork bifurcation for the Coanda effect via point-wise evaluation of the vertical velocity in the channel's centerline, and the corresponding coexisting solutions for the symmetric and asymmetric branches at $Re = 156$, i.e.\ $\mu = 0.5$.}
    \label{fig:pitchfork_coanda}
\end{figure}

Let us consider the steady and incompressible Navier-Stokes equations given by
\begin{equation}
\label{eq:NS_eq}
\begin{cases}
-\mu \Delta v + v\cdot\nabla v + \nabla p=0, \quad &\text{in} \ \Omega, \\
\nabla \cdot v = 0, \quad &\text{in} \ \Omega, \\
v = v_{\text{in}}, \quad &\text{on} \ \Gamma_{\text{in}}, \\
v = 0, \quad &\text{on} \ \Gamma_{\text{w}}, \\
- pn + (\mu \nabla v) n = 0, \quad &\text{on} \ \Gamma_{\text{out}},

\end{cases}
\end{equation}
where $v$ and $p$ represent the velocity and the pressure of the fluid, respectively, and $\mu$ is the kinematic viscosity. Let $\Gamma_{\text{in}} = \{0\}\times[2.5, 5]$ be the portion of the boundary $\partial \Omega$ where non-homogeneous Dirichlet conditions are imposed, namely via the inlet velocity $v_{\text{in}} = [20(5-x_2)(x_2 -2.5), 0]^T$,  $\Gamma_{\text{out}} = \{50\}\times[0, 7.5]$ be the ``free-flow" boundary at the outlet, and $\Gamma_{\text{w}} = \partial \Omega \setminus \{\Gamma_{\text{in}} \cup \Gamma_{\text{out}}\}$ the wall region featuring homogeneous Dirichlet boundary conditions.
The flow regime is characterized by the dimensionless Reynolds number defined as $\text{Re} = Uh / \mu$, i.e., the ratio between inertial and viscous forces, where $U$ and $h$ are the characteristic velocity of the flow and the characteristic length of the domain $\Omega$. In this specific case, the parameters are given by: $U = 31.25$ (the maximum inlet velocity) and $h = 2.5$ (the channel section). 
Investigating the flow behavior for different values of $\mu$, i.e., spanning different flow regimes, the model exhibits a \emph{pitchfork bifurcation} \cite{seydel2009practical}. Within this setting, it is well-known that the system becomes ill-posed at the bifurcation point located at $\mu^* \approx 0.96$ \cite{PichiArtificialNeuralNetwork2023,PintoreEfficientComputationBifurcation2021,khamlich2021model,BravoGeometricallyParametrisedReduced2024,TonicelloNonintrusiveReducedOrder2024}. In particular, the flow goes from exhibiting a unique symmetric state to admit multiple coexisting ones for the same value of the parameter $\overline{\mu} \leq \mu^*$, while for $\mu > \mu^*$ a unique laminar symmetric profile occurs \cite{QuainiSymmetryBreakingPreliminary2016,PichiDrivingBifurcatingParametrized2022a}. In fact, below such a critical threshold, the flow tends to attach either to the upper or lower wall boundary, breaking the symmetry of the configuration. Moreover, the symmetric behavior still exists in the bifurcating regime, but exchanges its stability properties with the asymmetric branches. That is, for the same value of the parameter $\overline{\mu}< \mu^*$, three solutions coexist. This phenomenon is represented in Figure \ref{fig:pitchfork_coanda}, where a bifurcation diagram is depicted for the considered test case, together with the three coexisting solutions for $\mu=0.5$.
To investigate the loss of uniqueness in a neighborhood of the pitchfork bifurcation, we set the parameter space as $\mathcal{P} = [0.5, 2.0]$, inducing a parametrization on the Reynolds number that varies in the range $Re\in [39, 156]$. 

Let us define 
$\V=\left(H^1(\Omega)\right)^2$, $\V_{\text{in}}=\{v \in \V \mid v=v_{\text{in}} \text{ on }\Gamma_{\text{in}}, v=0 \text{ on }\Gamma_{\text{wall}}\}$, $\V_0=\{v \in \V \mid v=0 \text{ on }\Gamma_{\text{in}} \cup \Gamma_{\text{wall}}\}$ and $\Q=L^2(\Omega)$. The weak formulation of \eqref{eq:NS_eq} reads: for a given $\mu \in \mathcal{P}$, find $v \in \V_{\text{in}}$ and $p \in \Q$ such that
\begin{equation}
\label{eq:gal_ns2}
\begin{cases}
a(v,\psi; \mu) +c(v,v,\psi) +b(\psi,p) = 0, \quad &\forall \, \psi \in \V_0, \\
b(v,\pi) = 0, \quad &\forall \, \pi \in \Q ,
\end{cases}
\end{equation}
where the weak formulation terms are defined for all $v$, $w$, $\psi \in \V$ and $p \in \Q$ as
\begin{equation*}
\label{eq:forms}
a(v, \psi; \mu) =\mu\int_\Omega\nabla v\cdot\nabla \psi \, d\Omega, \qquad
b(v, p) = -\int_\Omega(\nabla\cdot v) \hspace{.05cm}p \, d\Omega,
\end{equation*}
and 
$$
c(v, w, \psi)=\int_\Omega \left(v\cdot\nabla w\right) \psi \, d\Omega. 
$$
For the numerical discretization, we consider a mesh on the domain $\Omega$ and exploit the FE method with Taylor-Hood pair $\mathbb{P}^2$-$\mathbb{P}^1$, resulting in $N_{h} = N_h^v + N_h^p = 24301$ degrees of freedom.

At the reduced level, we build the reduced spaces $\mathbb V_{N_v}$ and $\mathbb Q_{N_p}$ for velocity and pressure, respectively of dimensions $N_v$ and $N_p$, by means of the developed greedy strategies, project the system on these reduced coordinates and solve the reduced, possibly deflated, problem via the Newton method. 

To guarantee the well-posedness of the Navier-Stokes system at the reduced level we enrich the velocity space by means of the \emph{supremizer stabilization}, see \cite{rozza2007stability} for details. This approach enlarges the velocity space to satisfy the reduced inf-sup condition \eqref{eq:inf_sup_N}, so that ${N_v} = 2N$ and ${N_p} = N$, where $N$ is the number of basis functions for each unknown.

\subsection{The adaptive-greedy performances} One of the main difficulties when dealing with bifurcating problems is that 
to capture the sudden change in the solution's behavior, one is usually restricted to a fine discretization of the parametric space. This also forces the HF sampling to waste computational time in the uniqueness region where the system could be easily approximated via only a few sampling points, instead. To perform an efficient and unbiased discovery of the parametric dependency, i.e., with no prior information concerning the location of the bifurcation point, we test here the performance when coupling the vanilla-greedy algorithm with the adaptive strategy\footnote{We remark that such adaptive strategy is not restricted to the greedy algorithm and may be coupled to other reduction approaches.} developed in Section \ref{sec:adaptivegreedy}. In this experiment,
we chose a maximum number of basis functions $N_{\text{max}}=35$, and a tolerance $\varepsilon = 10^{-3}$. To efficiently discover the bifurcating location via the reduced inf-sup indicator, we started from a very coarse discretization of the parametric space with only four equispaced values in the parametric interval $\mathcal{P}_h$. Once enough information is captured by the reduced model, i.e., the extrema of $\mathcal{P}_h$ are represented in the basis, we perform an equispaced refinement with $n_{\text{ref}}=4$ points in a neighborhood of the approximation of the bifurcation point given by the
parameter minimizing 
the stability indicator, i.e.\ $\mu_{\text{bif}} = \argmin_{\mu \in \mathcal{P}_h} \beta_{N}^h(\mu)$ defined in Section \ref{sec:adaptivegreedy}. 

\begin{figure}[h]
    \centering
    \includegraphics[width=0.8
\textwidth]{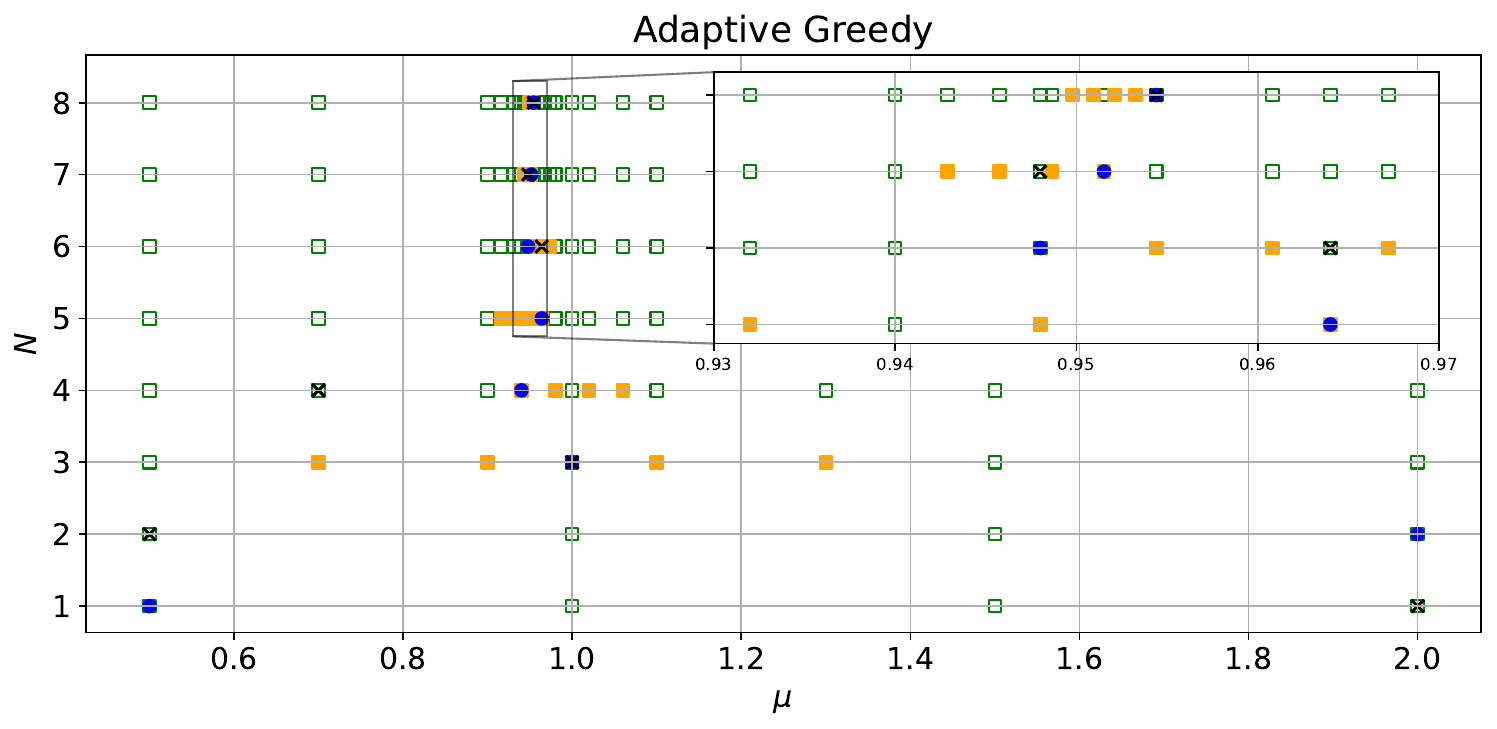}
    \caption{Evolution of the refinement strategy for the parametric space versus the number of basis/iterations.}
    \label{fig:mugreedy}
\end{figure}

The stopping criterion for the adaptive strategy is 
incremental 
with tolerance fixed as $\mathsf{tol} = 10^{-2}$. 
In Figure \ref{fig:mugreedy}, we show the evolution of the training set as the number of iterations of the greedy algorithm increases. Green squares denote the initial and retained training points from the previous iteration, while the filled yellow squares represent the equispaced points we add at each iteration.
Moreover, we indicate with the blue circle the current approximation of $\mu_{\text{bif}}$
around which the refinement is performed. Finally, the black crosses represent the parameter chosen by the greedy procedure, for which the corresponding snapshot is computed and added to the reduced basis. As anticipated, starting from a very sparse dataset, we considerably reduce the computational complexity of the offline phase, and allow for an efficient investigation of the parametric space. At the same time, the adaptive strategy concentrates the samples around the bifurcation point $\mu^* \approx 0.96$. It is worth noticing that such a procedure is very informative for the greedy strategy since it refines the sampling where the solution manifold exhibits a non-smoothness behavior and needs a finer discretization to accurately approximate the bifurcating phenomenon.

\begin{figure}[h]
    \centering
    \includegraphics[width=0.4\textwidth]{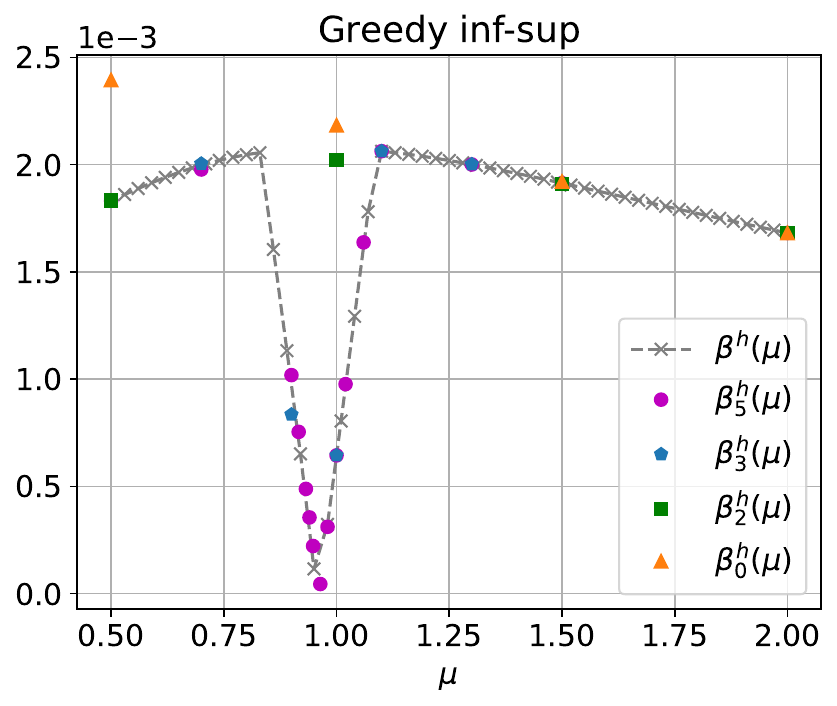}
    \caption{$\beta_N^{h}(\mu)$ approximation w.r.t.\ different iterations for the adaptive-greedy strategy compared with $\beta^h(\mu)$ computed on a fine equispaced discretization of $\mathcal{P}$.}
    \label{fig:betacomp}
\end{figure}

Figure \ref{fig:betacomp} shows the approximation of $\beta_N^{h}(\mu)$ compared to $\beta^{h}(\mu)$ computed on a fine equispaced discretization of $\mathcal{P}$ for several iterations of the adaptive-greedy strategy.
The refinement produces an automatic and efficient sampling that concentrates the points around $\mu^*$, thanks to the information coming from $\mu_{\text{bif}}$, adaptively matching the HF inf-sup constant with the reduced one. 
As specified in Section \ref{sec:adaptivegreedy}, this procedure can be exploited in conjunction with a deflation strategy. In particular, it could be considered as a pre-processing phase prior to applying the deflation technique to the system only for $\mu < \mu^*$, alleviating the computational costs of the deflated-greedy procedure. However, given the complexity of the deflation strategy itself, in the following, we discuss its properties without exploiting the adaptive approach only as a preliminary step, while focusing on the approximation capabilities of the deflation method both at the full and reduced order levels.


\subsection{The deflated-greedy performances}
In this section, we show the capability of the deflated-greedy strategy to certify multiple branches while building a reduced order model for a comprehensive approximation of the bifurcating phenomenon. The Coanda effect, in the parametric regime considered, features three coexisting solutions, meaning the number of branches is $K=3$. For the algorithm, we chose a maximum number of basis functions $N_{\text{max}}=35$, while the training set is discretized via an equispaced sampling with $|\mathcal P_h|=51$. The deflated-greedy tolerance is set to $\varepsilon=10^{-3}$, and the initial parameter for the strategy is taken as $\mu_0 = 2$. 
In particular, we employ as deflation parameters $r=2$ and $\sigma = 1$, as it is a common choice in the literature, but different bifurcating PDEs may need ad-hoc values \cite{Farrell2015A2026}. 
As previously discussed, a nice property of the deflation technique is that it does not require prior assumptions to set the initial guess for the nonlinear solver. Thus, both Newton methods, HF and RB, are initialized via the trivial null guess. Moreover, the stopping criterion for the standard and deflated Newton methods is imposed by setting a tolerance $\mathsf{tol}=10^{-10}$ over the residual vectors. 

To provide a comprehensive comparison between different reduction approaches, as a measure of performance we consider the relative reduced error defined as
\begin{equation}
    \label{eq:error_rel} E(\mu) = E_v(\mu) + E_p(\mu) = \frac{\norm{v_h(\mu) - v_N(\mu)}_{\mathbb V}}{\norm{v_h(\mu)}_{\mathbb V}} + \frac{\norm{p_h(\mu) - p_N(\mu)}_{\mathbb Q}}{\norm{p_h(\mu)}_{\mathbb Q}}, 
\end{equation}
and computed on the test set $\mathcal P_\text{te}$
of equispaced points in $\mathcal P$ of cardinality $|P_\text{te}| = 151$.

\begin{figure}[htbp]
    \centering

    \caption*{\underline{POD}}
    \begin{subfigure}[h]{0.3\textwidth}
\includegraphics[width=\textwidth]{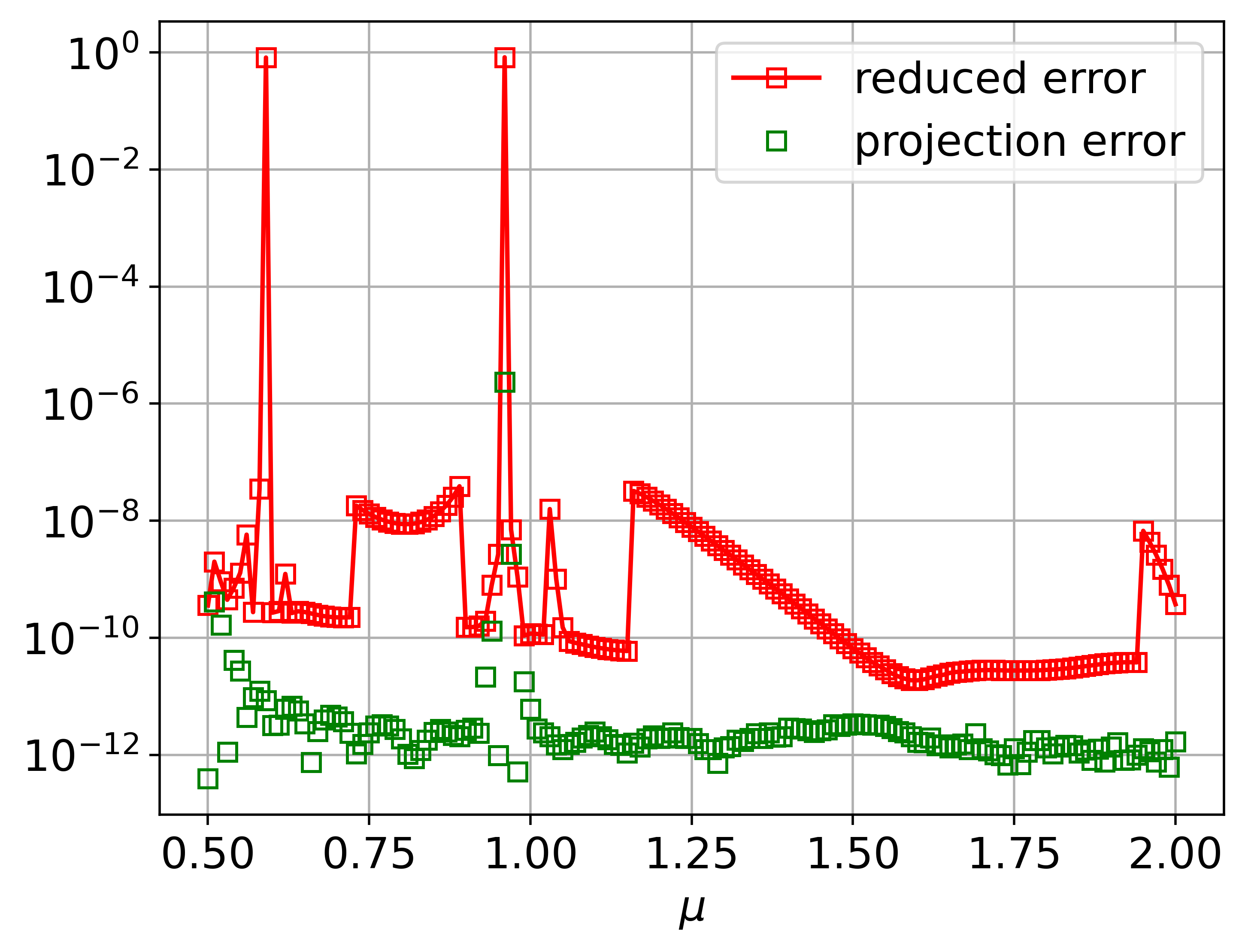}
        \caption{sym.}
        \label{sfig:a}
    \end{subfigure}
    \begin{subfigure}[h]{0.3\textwidth}\includegraphics[width=\textwidth]{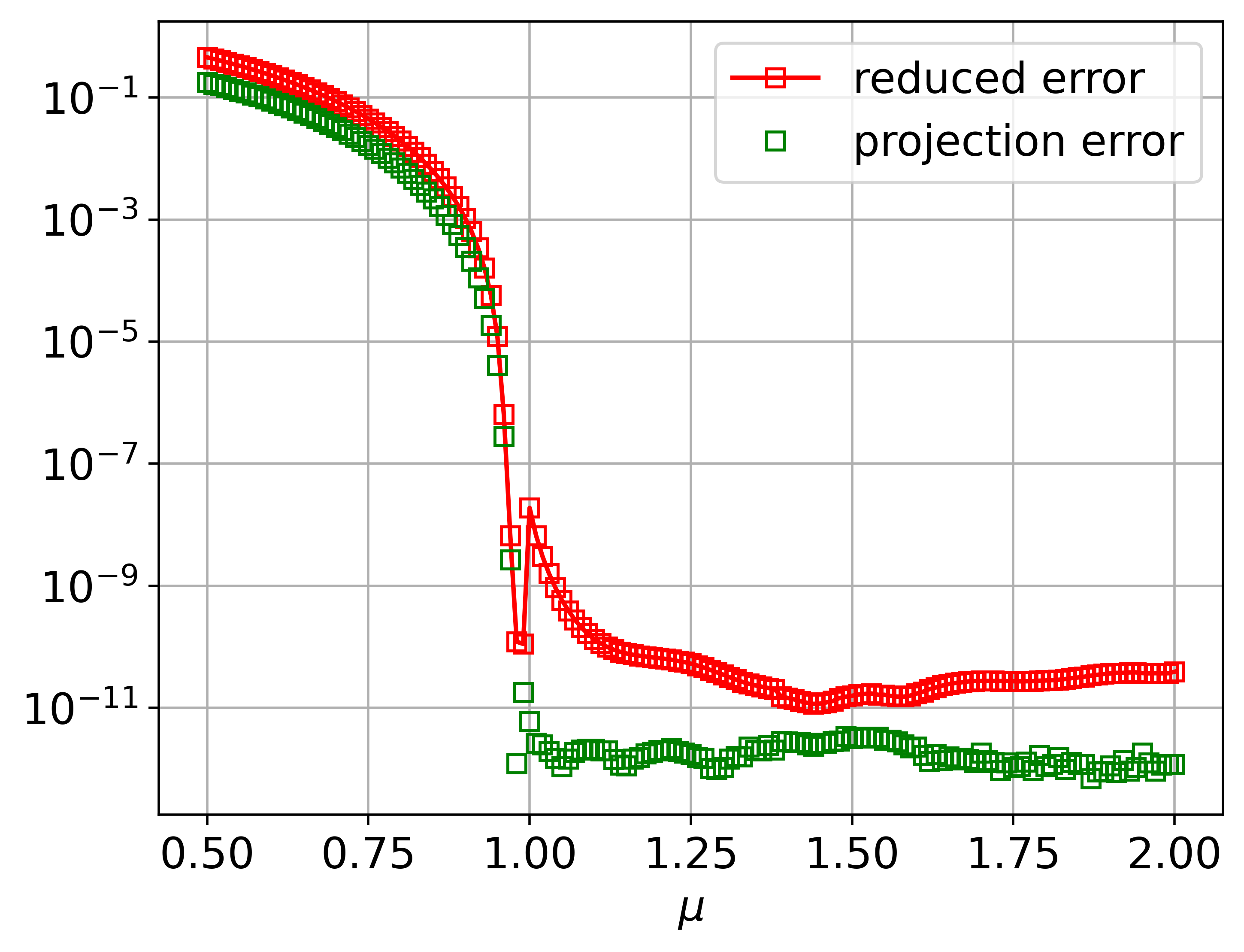}
        \caption{asym.\ upper}
    \end{subfigure}
    \begin{subfigure}[h]{0.3\textwidth}\includegraphics[width=\textwidth]{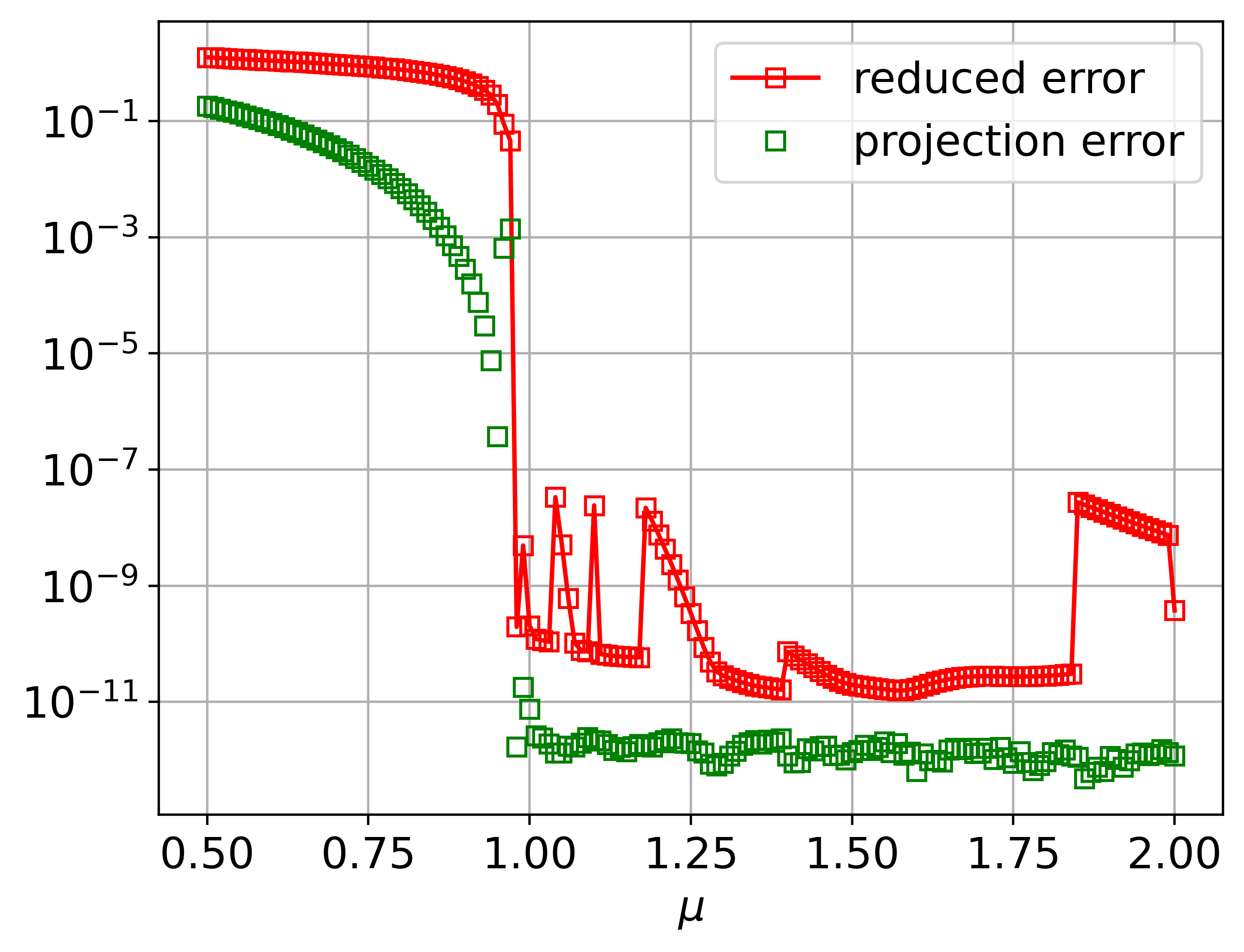}
        \caption{asym.\ lower}
    \end{subfigure}

        \caption*{\underline{vanilla-greedy}}
    \begin{subfigure}[h]{0.3\textwidth}
\includegraphics[width=\textwidth]{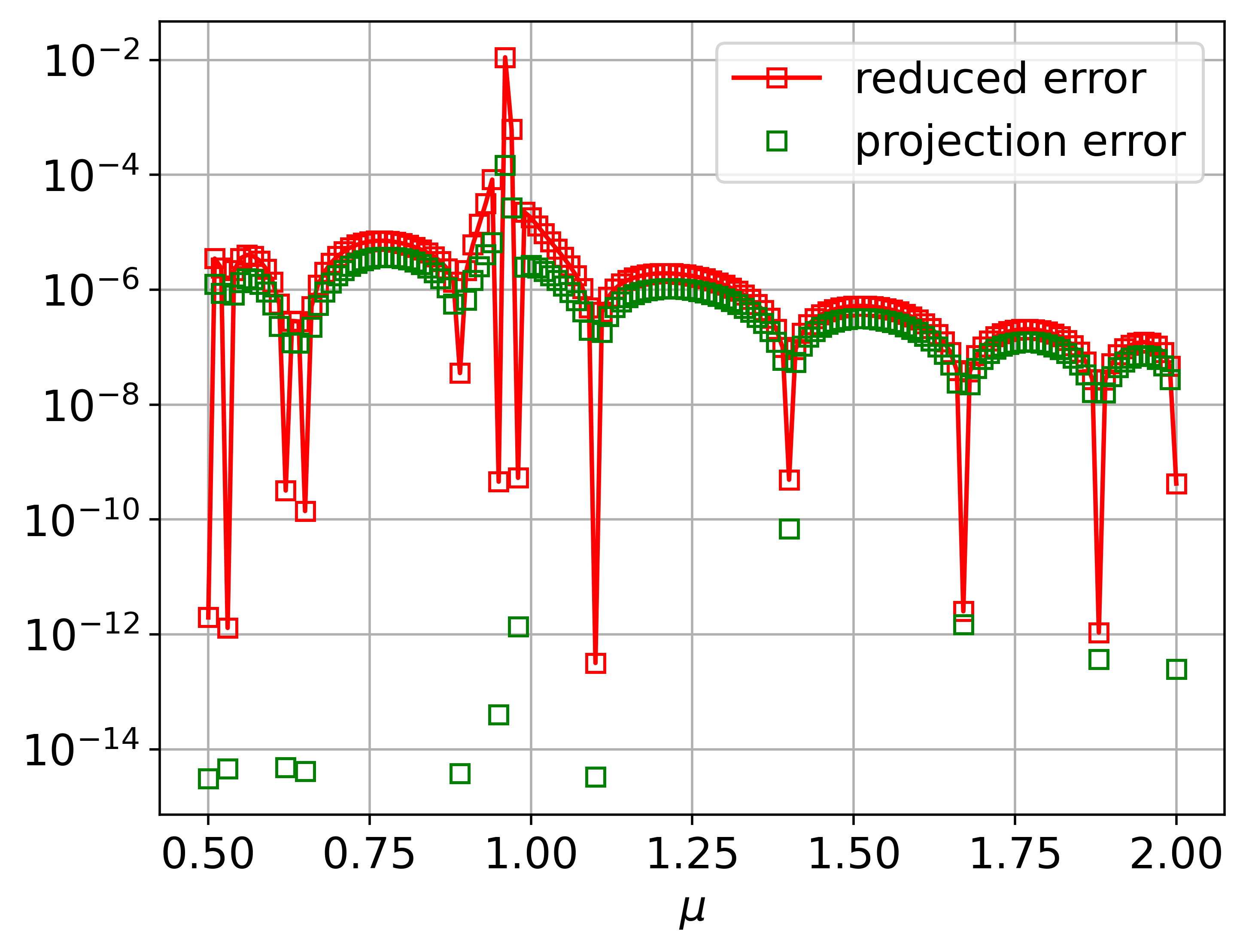}
        \caption{sym.}
    \end{subfigure}
    \begin{subfigure}[h]{0.3\textwidth}\includegraphics[width=\textwidth]{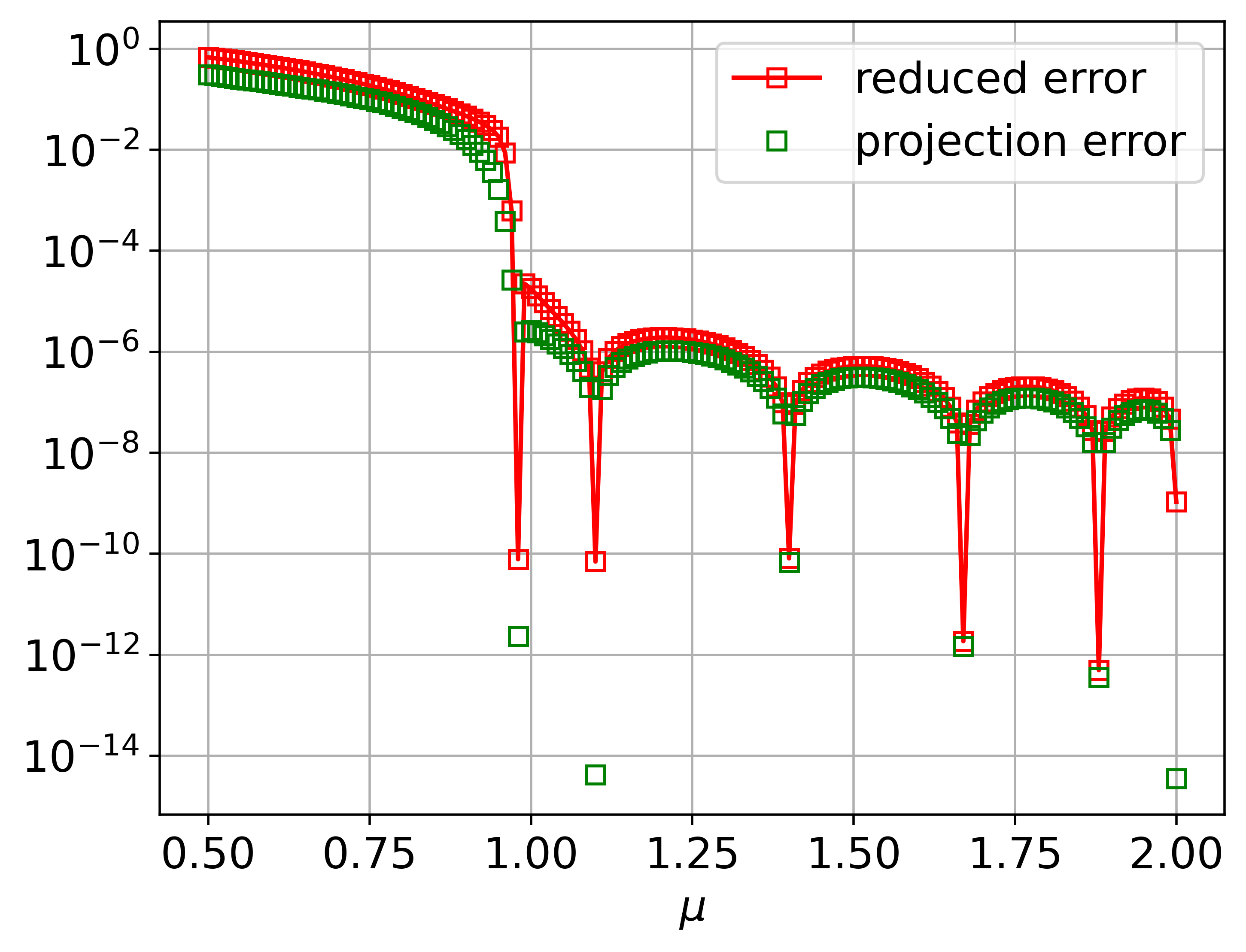}
        \caption{asym.\ upper}
    \end{subfigure}
    \begin{subfigure}[h]{0.3\textwidth}\includegraphics[width=\textwidth]{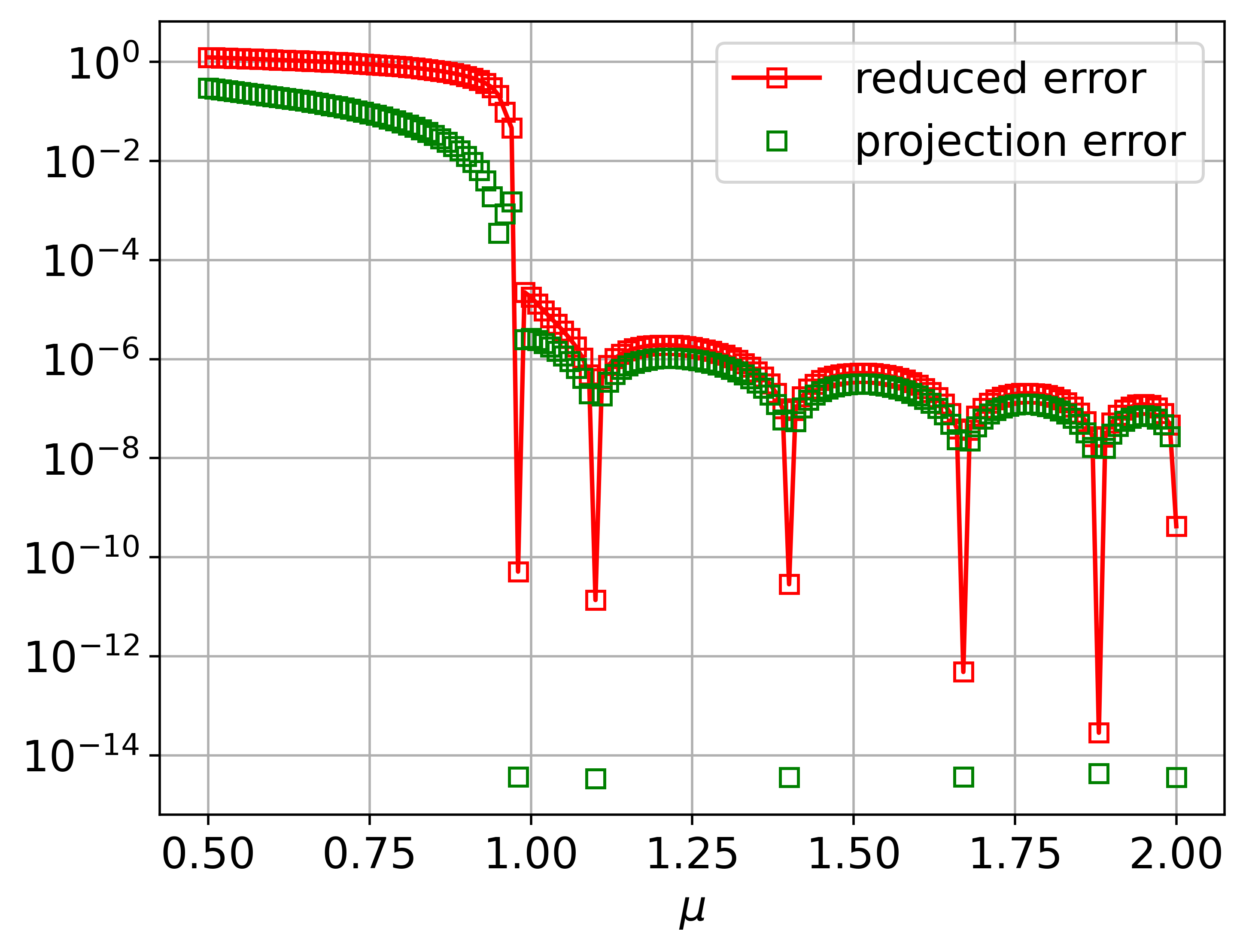}
        \caption{asym.\ lower}
    \end{subfigure}

    \caption*{\underline{deflated-greedy}}
    \begin{subfigure}[h]{0.3\textwidth}
\includegraphics[width=\textwidth]{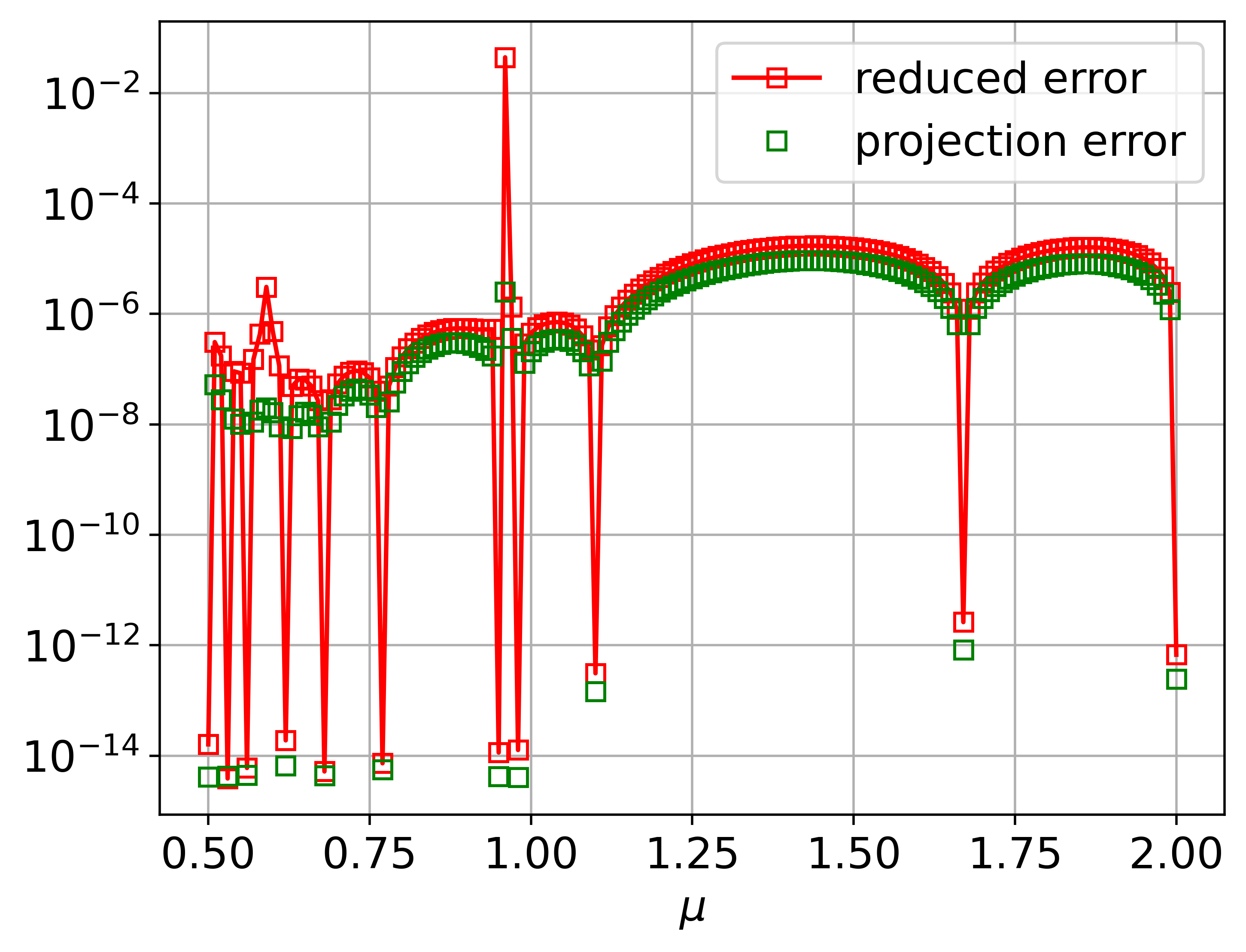}
        \caption{sym.}
        \label{sfig:g}
    \end{subfigure}
    \begin{subfigure}[h]{0.3\textwidth}\includegraphics[width=\textwidth]{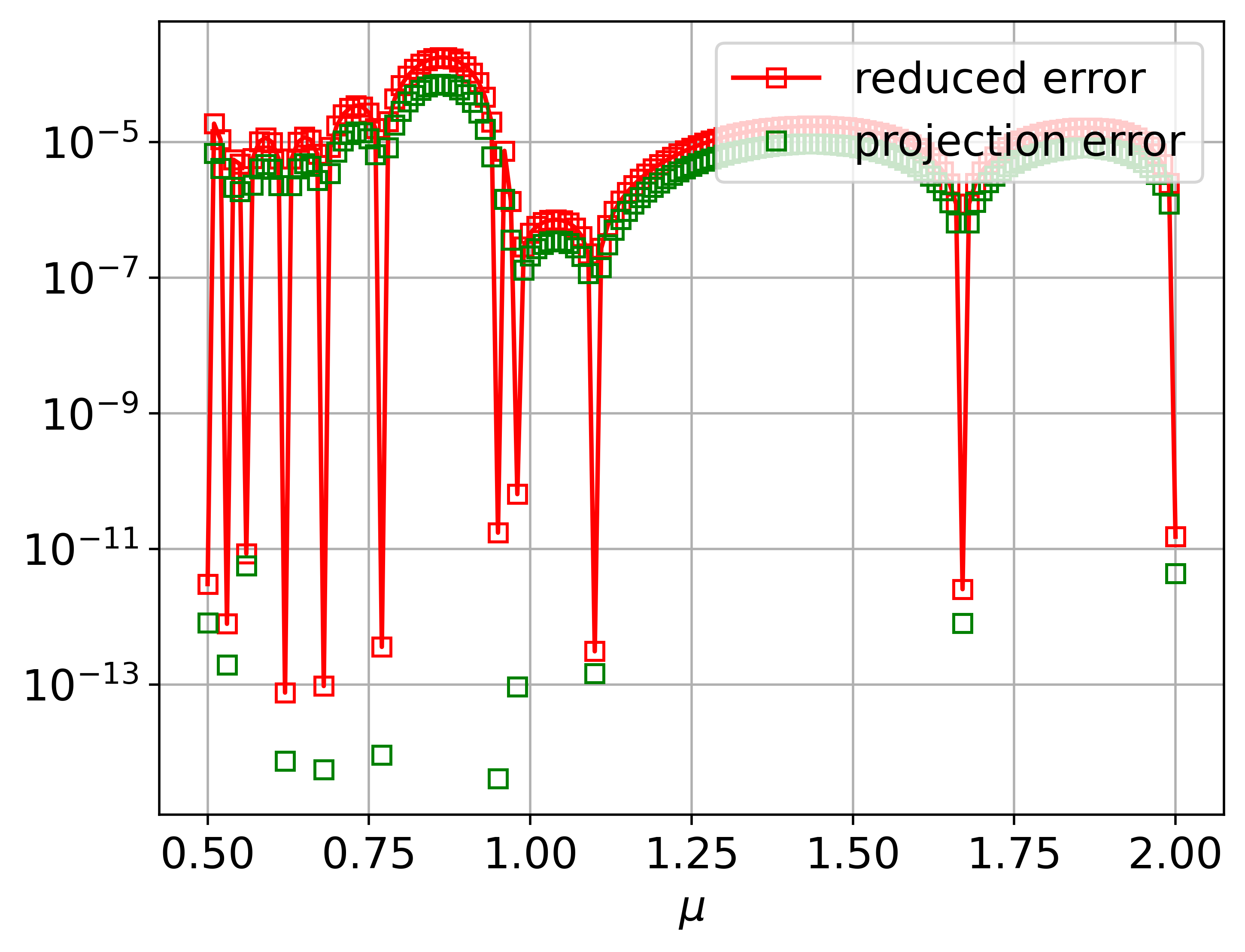}
        \caption{asym.\ upper}
    \end{subfigure}
    \begin{subfigure}[h]{0.3\textwidth}\includegraphics[width=\textwidth]{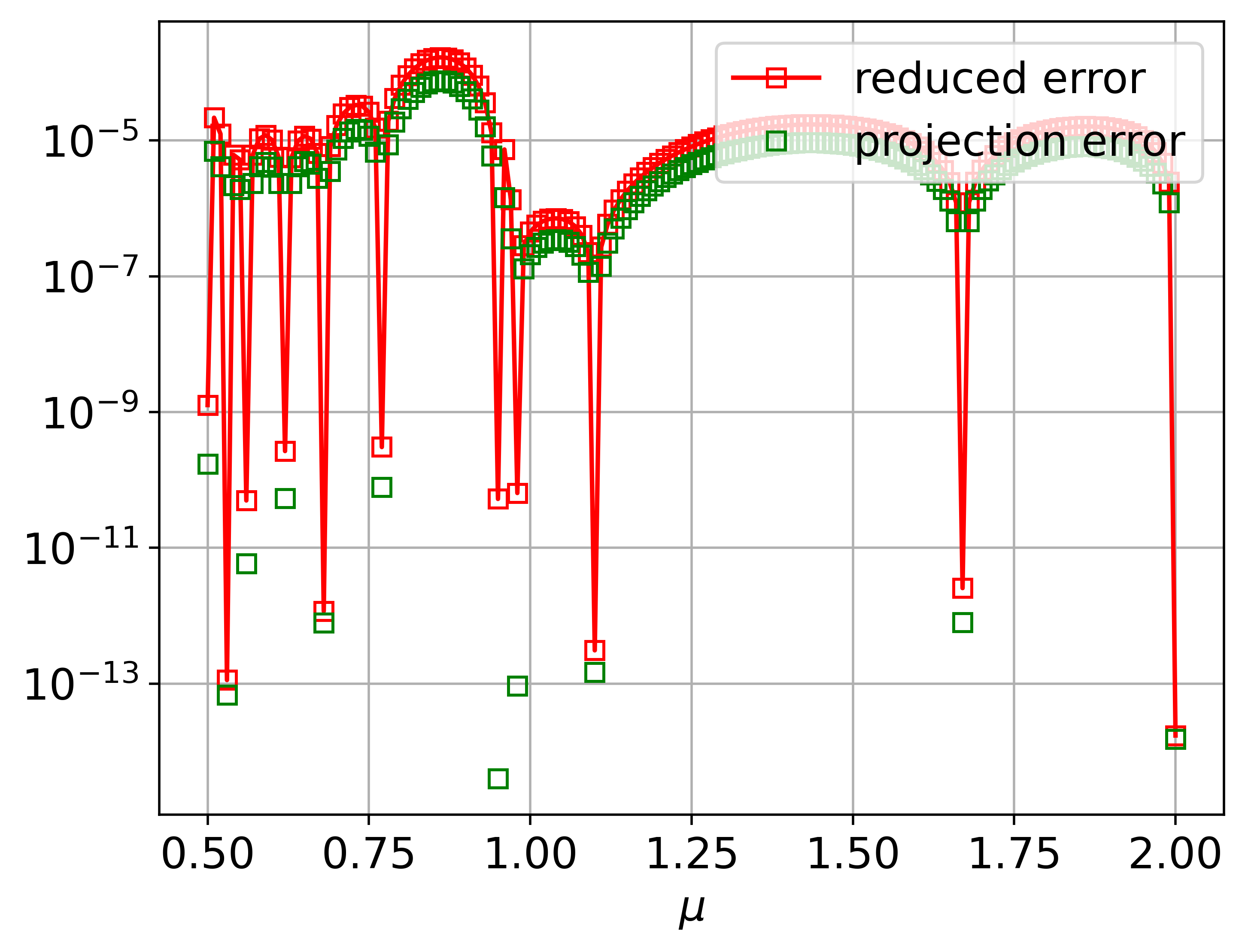}
        \caption{asym.\ lower}
    \end{subfigure}

    \caption{Relative error and  branch certification. Symmetric, upper asymmetric, and lower asymmetric solutions from left to right. POD ($N=25$), vanilla-greedy ($N=12$), and deflated-greedy ($N=25$), are depicted respectively in the top, center, and bottom rows.}
    \label{fig:certification_2}
\end{figure}

We now show in Figure  \ref{fig:certification_2} the accuracy and certification properties of POD, vanilla-greedy, and the proposed deflated-greedy approach for the simultaneous reduced reconstruction of the pitchfork bifurcation comprising the symmetric and two asymmetric branches.

Within the configuration setting described above, to compute the performance in the parametric space, we retained the first $N=25$ eigenvectors for the POD data compression, in order to have a fair comparison with the deflated-greedy strategy, for which the stopping criterion on the estimator set is also met for $N=25$. In contrast, the results for the vanilla-greedy method have been obtained with $N=12$, since the convergence criterion, with $\varepsilon = 10^{-5}$, has been reached much in advance for the sampled data. While for general problems such behavior would be favorable, in such a context is already a clear symptom that the vanilla-greedy strategy is completely missing the coexistence of the states, not retrieving the full information from the system.

In particular, in Figure \ref{fig:certification_2} we depict the relative error $E(\mu)$ obtained both via direct projection of the solutions on the constructed basis (green squares), and via the standard Galerkin projection (red squares) solving the reduced nonlinear system.

From the first row, showing the accuracy of the branch-wise POD (based on the solutions obtained from the unbiased Newton solver) in approximating the three branches, it is clear that, even though the symmetric configurations have been captured almost perfectly by the reduced basis, the same is not true for the two asymmetric ones, for which the error grows exponentially fast in the non-uniqueness regime, and in the direction of lower viscosity values. The agreement of the projection and reduced errors confirms the impossibility of POD to effectively capture the bifurcating phenomenon at the reduced level (if no additional assumptions are made about its existence). Indeed, since the projection error represents the best approximation in the reduced space spanned by the basis functions, the possible differences among the two performance indicators are only due to the Galerkin step or to the convergence of the reduced Newton convergence history. This is, for example, the case of Figure \ref{sfig:a} where the reduced Netwon technique did not converge near the bifurcation point and for a high Reynolds number, but the errors for projected solutions are still always below $10^{-5}$, denoting that the reduced basis is correctly approximating (only) the symmetric branch. 

When applying the vanilla-greedy approach, the sampling in $\mathcal{P}_h$ is driven by the error estimator, guiding the model toward the discovery of less represented parametric regions. The main issue in the bifurcating context is that a portion of the parametric space could be well described by the available basis but only for a subset of the coexisting branches, preventing the error certification.

From the second row in Figure \ref{fig:certification_2}, we observe the usual error evolution in the parametric space for greedy techniques relative to the direct sampling of the basis as the snapshots of the system.
Once again, the symmetric branch has been recovered, but the bifurcating regime has been missed.
Indeed, 
even if the vanilla-greedy approach reaches the stopping criterion already for $N=12$, the standard estimator is computed through the reduced representations of the approximated branches, and in this case, the bifurcating solutions have not been sampled. 
Moreover, even the certification of the symmetric branch is not guaranteed due to the inherently ill-posed nature of the model around the critical location.

Having detailed the issues of POD and vanilla-greedy approaches, we finally comment on the results of the certified approximation obtained via the proposed methodology. As it can be seen in the last row of Figure \ref{fig:certification_2}, for all branches and all parameters
in $\mathcal{P}_\text{te}$, the reduced error is reliably below the tolerance bound $\varepsilon=10^{-3}$, meaning that only the deflated-greedy strategy is capable of certifying all the admissible branches coexisting in the parametric range considered.

The unique point in the parametric space that seems to violate the certification is indeed the one corresponding to the bifurcation point while reconstructing the symmetric branch in Figure \ref{sfig:g}, and the reason is that, as before, the reduced Newton algorithm converged to a different branch. Thus, as confirmed by the projection error, the issue is not related to the basis, but to the sensitivity of the nonlinear solver in a neighborhood of the critical point.

Moreover, we highlight that, in addition to the increased approximation capability and the certification of the reduced model, the deflated-greedy approach inherently benefits from the parsimonious sampling of the high-fidelity manifold. Contrarily to POD approaches, which require a biased and complete dataset (i.e., with the knowledge of the bifurcating phenomenon), and vanilla-greedy, 
the deflated reduced approach 
samples more frequently in the bifurcating regime where the branches coexist, and thus where more information is needed.

These observations fully comply with the features of the proposed algorithm: paying the price of the deflation step, the methodology can detect potential non-uniqueness behavior without prior knowledge, exploiting the greedy estimators to drive the sampling, providing more information on the bifurcation, and certifying the discovered branches by adding the corresponding basis functions to the reduced space. 

\begin{figure}[h]
    \centering
    \includegraphics[width=\textwidth]{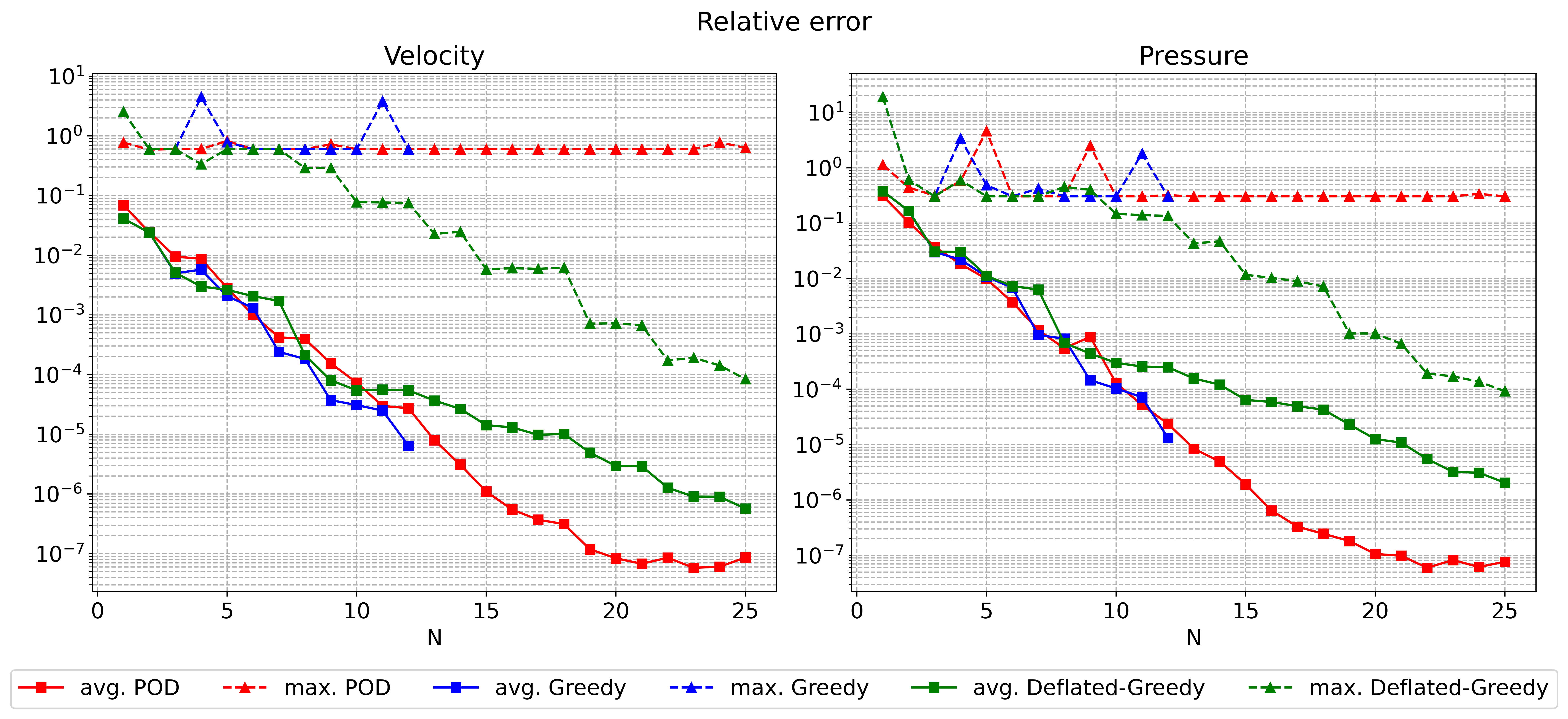}
    \caption{Average and maximum relative error comparison w.r.t.\ $N$ between POD (red lines), vanilla-greedy (blue lines) and deflated-greedy (green lines) algorithms for the velocity and pressure asymmetric fields.}
    \label{fig:comparison}
\end{figure}

Finally, we remark that, as it has already been noticed when applying POD approaches to different bifurcating phenomena \cite{pichi_phd,PichiDrivingBifurcatingParametrized2022a}, the peaks of the relative errors are usually located in correspondence to the critical regions, providing an \emph{a-posteriori} indicator of the bifurcation points. 
In Figure \ref{fig:comparison}, we compare the average and maximum relative errors $E_v$ and $E_p$  w.r.t.\ the dimensionality of the reduced space for deflated-greedy, vanilla-greedy, and branch-wise POD. We show the results for the asymmetric lower solution, i.e., the bifurcating branch  missed by both POD and vanilla-greedy, since the symmetric behavior is well approximated by all the strategies. We see that for the three approaches the average relative errors decrease w.r.t.\ $N$. POD and vanilla-greedy even seem to perform better than the deflated-greedy in terms of average error for $N=12$, but we remark that here the two approaches contain only the information from the symmetric branch, representing almost exactly the uniqueness behavior, consisting in almost 2/3 of the entire parametric range. Thus, they are more expressive for a fixed amount of basis functions than the deflated-greedy, which instead contains the knowledge of all discovered branches. 
However, the error certification features and the superior performance of the deflation-augmented technique are clear when comparing the behavior of the maximum error and increasing the number of basis functions. In fact, both the POD and the vanilla-greedy approaches fail to certify the model for each parameter, while the deflated-greedy strategy shows a comparable trend for both the maximum and average errors, eventually reaching the stopping certification criterion and providing a reliable model. 
Such a nice performance comes from the deflation step, which allows us to discover and add qualitatively different basis functions corresponding to different coexisting solutions for the same parametric instance. Furthermore, as concerns the POD approximation, due to the lack of certification, a clear saturation effect in the average error is observed for $N \geq 20$.

To further clarify the main differences between the two greedy approaches, especially in terms of the error bound, we investigate the behavior of the \textit{vanilla} estimator $\Delta_N(\mu)$ and of the \textit{deflated} estimator $\max{\Delta}$.

\begin{figure}[H]
    \centering
    \includegraphics[width=0.49\textwidth]{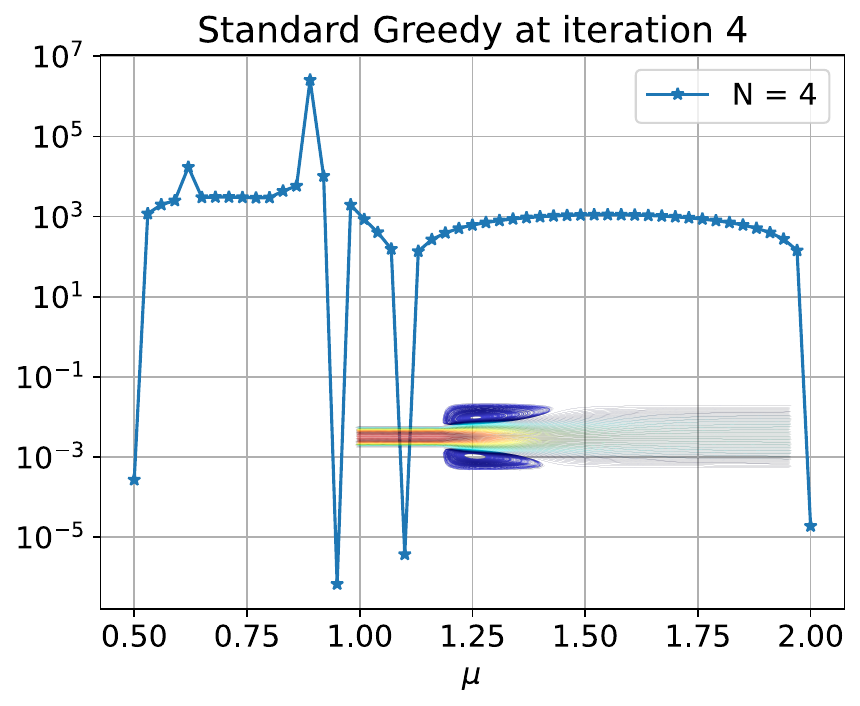}
    \includegraphics[width=0.49\textwidth]{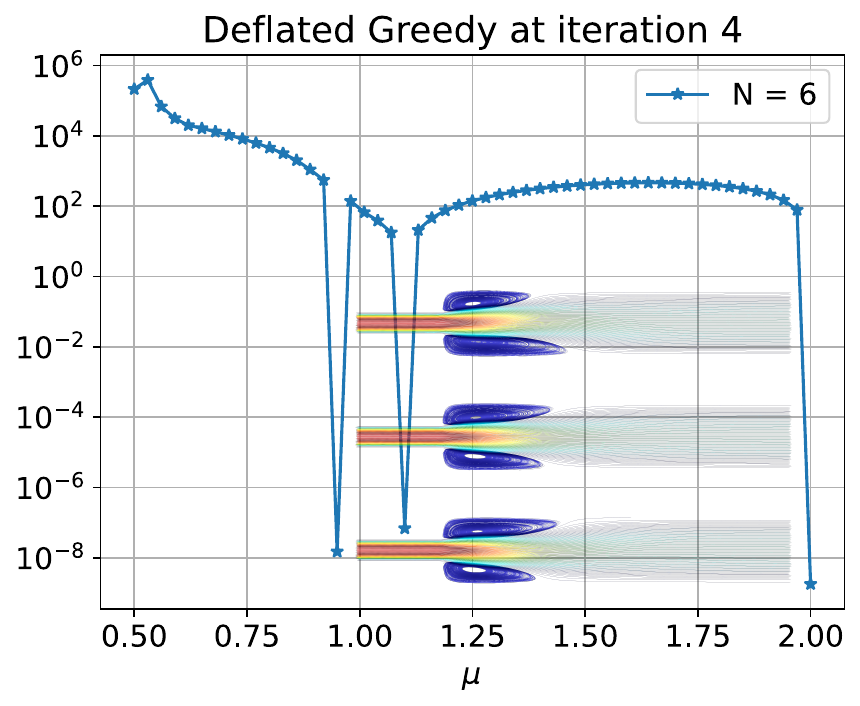}
    
    \caption{Bound for deflated-greedy and vanilla-greedy algorithms at the fourth iteration. The plot depicts the sampled snapshot(s) and the effect on the detection of different admissible branches.}
    \label{fig:boundandsnap}
\end{figure}
In Figure \ref{fig:boundandsnap}, we compare the fourth iteration of the deflated- and vanilla-greedy algorithms, i.e., the first one for which the two approaches behave differently. \B{On the abscissa axis, we represent the parameter value, and on the ordinate axis, we provide the value of the linear error estimator bound for the vanilla-greedy and the maximum of the deflated estimator over the three branches.
The plot is enriched by the velocity streamlines of the sampled snapshots for the two procedures: one for the vanilla-greedy and three for the deflated-greedy.}
In the left plot, corresponding to the vanilla-greedy, it can be seen that the greedy strategy selects $\mu_4 = 0.95$, the basis has dimension four, and the bound is minimized at the sampled locations. The seemingly ``nice" behavior of the error bound actually hides the difficulty of the approach to provide information on the bifurcating phenomenon. Indeed, as it can be observed from the streamline of the flow exhibiting the vortex pattern, even if $\mu_4 < \mu^*$, the sampled basis is symmetric, denoting that the whole reduced branch reconstruction, upon which the error bound is computed, only recovers the straight flux corresponding to the symmetric behavior.

In contrast, as depicted in the right plot, the deflated-greedy method at the same iteration is capable of discovering multiple coexisting branches through deflation. Thus, we add all possible information from the solution manifold, enrich the basis with $b_{\mu_4}=3$ snapshots corresponding to $\mu_4 = 0.95$, and encode both the symmetric and the two asymmetric configurations. Moreover, this has the desirable effect of detecting the reduced asymmetric branch. In fact, having sampled asymmetric profiles, and thus having the possibility to reconstruct the branching behavior, the maximum among the error bounds is no longer restricted to the representation of the symmetric branch and increases for lower viscosity values. This is indeed a confirmation that the wall-hugging behavior was not present in the basis, but being detected, it can now help the strategy to continue the asymmetric profiles towards the certification of the multiple branches.

\section{Conclusions}
\label{sec:conc}
In this work, we present two novel greedy algorithms for bifurcating nonlinear parametric PDEs in a ROM setting. 
We propose an adaptive-greedy strategy to detect the bifurcation point starting from sparse information on the parametric space.
Furthermore, we conceived a deflated-greedy
that certifies multiple branches employing the deflation strategy. Deflation is employed in two phases: (i) to enrich the reduced space with multiple snapshots with different physical behaviors, (ii) to evaluate multiple error estimators related to different reduced solutions. It guarantees a more informative reconstruction of the reduced solution and drives the research of the new snapshots toward the less-represented branch solution.
The two strategies have been tested in a sudden expansion channel flow, featuring a pitchfork bifurcation with three coexistent solutions. 
The adaptive-greedy can identify in a small number of iterations the bifurcation point with no previous knowledge of the system, while the deflated-greedy can certify all the branches. The results are compared to vanilla-greedy and POD in terms of accuracy w.r.t.\ the HF solutions.

Moreover, we investigated in 
\ref{app:estimator} the role and applicability of the nonlinear error estimator for different configurations of the Navier-Stokes problem.

This contribution is a first step towards an efficient and unbiased reduced order investigation of nonlinear bifurcating problems. The content of this work can be useful in many scientific and industrial fields, where strong nonlinear dynamics are of interest and need to be studied for several parametric instances without any knowledge of possible bifurcating behaviors. In particular, when considering tasks such as control, stability analysis, geometrical parametrization and multi-parameter design, that will constitute the main direction of future investigation.

\section*{Acknowledgments}
MS acknowledges the European Union's Horizon 2020 research and innovation program under the Marie Skłodowska-Curie Actions, grant agreement 872442 (ARIA). This study was carried out within the ``20227K44ME - Full and Reduced order modelling of coupled systems: focus on non-matching methods and automatic learning (FaReX)" project – funded by European Union – Next Generation EU  within the PRIN 2022 program (D.D. 104 - 02/02/2022 Ministero dell’Università e della Ricerca). This manuscript reflects only the authors’ views and opinions and the Ministry cannot be considered responsible for them. Moreover, MS thanks the INdAM-GNCS Projects ``Metodi di riduzione di modello ed approssimazioni di rango basso per problemi alto-dimensionali" (CUP E53C23001670001) and ``Metodi numerici efficienti per problemi accoppiati in sistemi complessi” (CUP E53C24001950001).
FP acknowledges the support provided by the European Union- NextGenerationEU, in the framework of the iNEST- Interconnected Nord-Est Innovation Ecosystem (iNEST ECS00000043 – CUP G93C22000610007) consortium and its CC5 Young Researchers initiative.


\bibliographystyle{abbrvurl}
\bibliography{main.bib}
\appendix
\begin{appendices}
\addtocontents{toc}{\protect\setcounter{tocdepth}{0}}
\section{A numerical comparison between estimators for bifurcating Navier-Stokes equations}
\label{app:estimator}
This Appendix focuses on the performances of the nonlinear \emph{a-posteriori} error estimator based on the Brezzi-Rappaz-Raviart theory \cite{BRR1} for the Navier-Stokes system, and the challenges arising in the bifurcating setting. 
\subsection{Definition of the error estimator} Let us provide here a more detailed explanation of the several terms appearing in the error estimator \eqref{eq:nonlinear_est} when dealing with the Navier-Stokes equations. In the considered framework, $K^h_N(\mu) = 2 \rho^2 M_c(\mu)$, where $\rho = \rho(\Omega)$ is the Sobolev embedding constant defined as
\begin{equation}
\rho^2 = \sup_{v \in \mathbb V} \frac{\norm{v}^2_{L^4(\Omega)}}{\norm{\nabla v}_{L^2(\Omega)}^2},
\end{equation}
and $M_c(\mu)$ is a function depending on the parametric problem, which in this case is given by $M_c(\mu)=\sqrt{2}$. To obtain an approximation of the Sobolev embedding constant $\rho$, one can exploit a strategy that combines an eigenvalue problem and a fixed-point algorithm, for further insight see \cite{Deparis20082039,ManzoniEfficientComputationalFramework2014}. 

An important quantity to look at in order to define the nonlinear estimator is $\tau_N(\bmu)$, which for the Navier-Stokes problem is defined as
\begin{equation}
    \label{eq:tau_man} 
    \tau_N(\bmu) = \frac{4\gamma(\rho; \bmu) \norm{\mathsf G( \mathsf B; \bmu)}_{\V_{} \times \Q_{}}}{\beta^h_N(\bmu)^2},
\end{equation}
where $\mathsf G$ is the algebraic residual of the system, $\gamma(\rho; \bmu) = \rho^2 M_c(\mu)$ is the continuity constant of the trilinear form $c(\cdot, \cdot, \cdot)$, and $\mathsf B$ is the global basis function matrix for velocity and pressure.
Thus, the specific nonlinear error estimator is given by
\begin{equation}
\label{eq:nonlinear_NS}
\Delta_N^{\text{nl}}(\mu) = \frac{\beta^h_N(\bmu)}{2 \gamma(\rho; \bmu) }\left (1 - \sqrt{1 - \tau_N(\bmu)} \right ),
\end{equation}
and it can be employed to certify the reduced error only when $\tau_N(\bmu) \leq 1$ for every $\mu \in \mathcal P_h$.
As anticipated in Section \ref{sec:estimator}, the condition $\tau_N(\bmu) \leq 1$ might be difficult to meet, especially for bifurcating phenomena. In fact, the existence of bifurcation points, at which the branching behavior takes place, causes the beta in-sup constant $\beta^h_N(\mu)$ to assume zero values at the critical points, and consequently large values of $\tau_N(\mu)$. 
In the following, we report some numerical results highlighting the comparison between the performances of the nonlinear error estimator \eqref{eq:nonlinear_NS} for the bifurcating sudden-expansion channel test a similar but non-bifurcating system. These results justify the exploitation of the linear bound \eqref{eq:linear_est} in the tests of Section \ref{sec:results}.

\subsection{Numerical results}
To further study the usage of the nonlinear error bound \eqref{eq:nonlinear_NS} in a bifurcating regime, we consider two numerical settings: (i) a slight modification of a flow over a backward-facing step proposed in \cite{ManzoniEfficientComputationalFramework2014} and (ii) the sudden-expansion channel test of Section \ref{sec:results}. As we have seen, in the parametric ranges considered, the latter admits coexisting solutions resulting from the bifurcating phenomena, while the former is a well-posed problem with a unique solution. For both settings, we employ the vanilla-greedy strategy comparing the behavior of $\tau_N(\mu)$ during the procedure.
In order to propose a fair comparison in terms of $\gamma(\rho; \mu)$ and $Re$ between the two test cases, we slightly modified the benchmark of the flow over a backward-facing step \cite{ManzoniEfficientComputationalFramework2014}. In particular, we adapt the geometry of (i) so that the two Sobolev embedding constants $\rho$ (which depend on the spatial domain) have comparable values, and we change the parameter space to also have comparable Reynolds numbers.

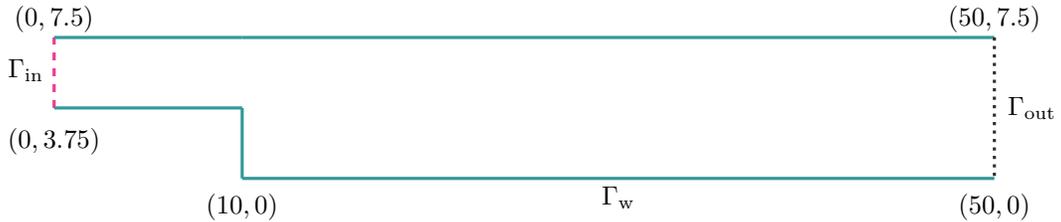
\begin{figure}[h]
\begin{center}
\begin{tikzpicture}[scale=.25]

\filldraw[color=magenta!90, very thick, dashed](0,3.75) -- (0.,7.5);
\filldraw[color=teal!80, fill=gray!10, very thick](0,3.75) -- (10,3.75);
\filldraw[color=teal!80, fill=gray!10, very thick](0,7.5) -- (10,7.5);
\filldraw[color=teal!80, fill=gray!10, very thick](10,3.75) -- (10,0);
\filldraw[color=teal!80, fill=gray!10, very thick](10,0) -- (50,0);
\filldraw[color=teal!80, fill=gray!10, very thick](10,7.5) -- (50,7.5);
\filldraw[color=black!80, fill=gray!10, very thick, dotted](50,7.5) -- (50,0);

\node at (-1.5,5.75){\color{black}{$\Gamma_{\text{in}}$}};
\node at (52,3.75){\color{black}{$\Gamma_{\text{out}}$}};
\node at (30,-1){\color{black}{$\Gamma_{\text{w}}$}};
\node at (0,2.){\color{black}{$(0,3.75)$}};
\node at (0,8.5){\color{black}{$(0,7.5)$}};
\node at (10, -1.5){\color{black}{$(10,0)$}};
\node at (50, -1.5){\color{black}{$(50,0)$}};
\node at (50, 8.5){\color{black}{$(50,7.5)$}};

\end{tikzpicture}
\end{center}
\caption{Domain $\Omega$. Dirichlet boundary conditions are applied over the dashed magenta line (non-homogeneous inlet) and the solid teal line (homogeneous). The dotted black line features 
{``free-flow"} boundary conditions.
}
\label{fig:domainManzoni}
\end{figure}

As concerns problem (i), we consider the Navier-Stokes equations on the domain $\Omega$ depicted in Figure \ref{fig:domainManzoni}.

Let $\Gamma_{\text{in}} = \{0\}\times[3.75, 7.5]$ and $\Gamma_{\text{out}} = \{50\}\times[0, 7.5]$ be the portions of the boundary $\partial \Omega$ where non-homogeneous Dirichlet and ``free-flow" boundary conditions are imposed, respectively. On the remaining boundary $\Gamma_{\text{w}} = \partial \Omega \setminus \{\Gamma_{\text{in}} \cup \Gamma_{\text{out}}\}$ homogeneous Dirichlet boundary conditions are applied.
The inlet velocity is $v_{\text{in}} = [21(7.5-x_2)(x_2 -3.75), 0]^T$. The parameter space is $\mathcal P = [0.024, 0.096]$, corresponding to $Re \in [39, 156]$. A $\mathbb P^2$-$\mathbb P^1$ {Taylor-Hood} approximation is used, with $N_h = N^v_h + N_p^h = 24386$. In this setting, the asymmetric domain $\Omega$, avoids recirculations and no bifurcation is obtained. Thus, we expect no issues using the nonlinear error estimator.
\begin{figure}[h]
    \centering
    \includegraphics[width=0.49\textwidth]{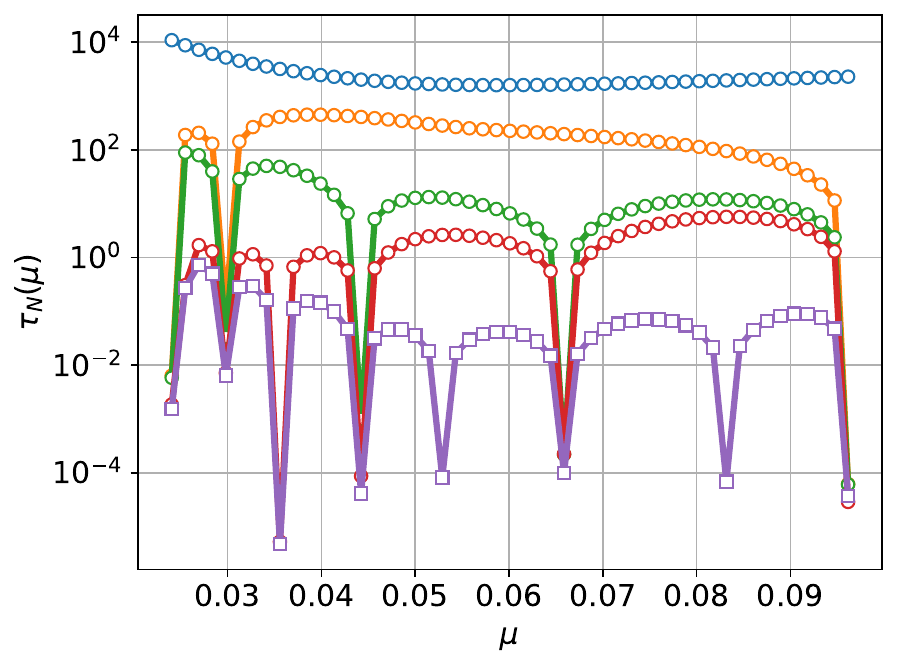}
    \includegraphics[width=0.49\textwidth]{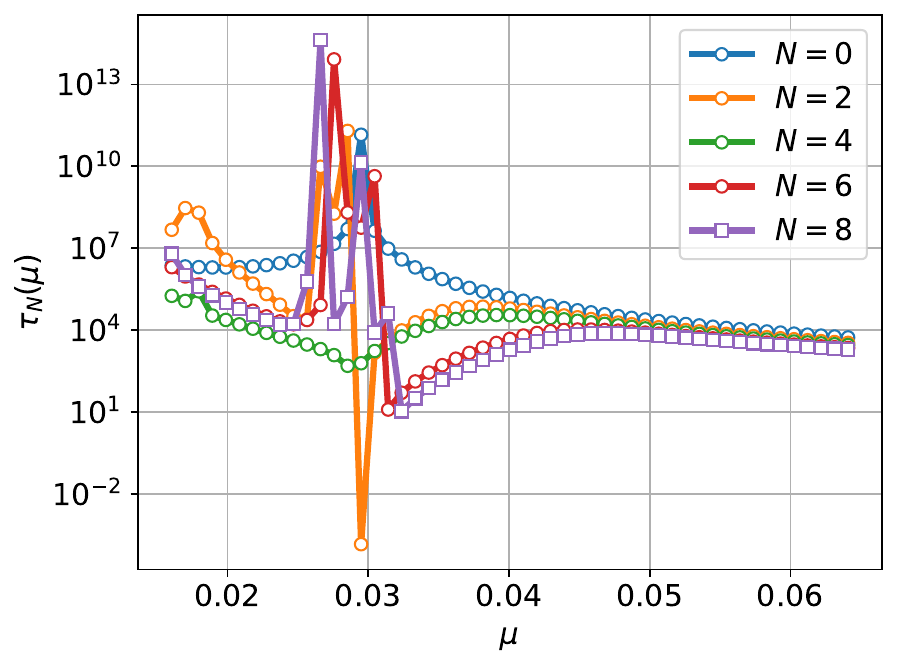}
    \caption{Evolution of $\tau_N(\mu)$ for several iterations of the vanilla-greedy for (i) the flow over a backward-facing step and (ii) the sudden-expansion channel, left and right plot, respectively.}
    \label{fig:Manzoni}
\end{figure}


We perform a vanilla-greedy procedure with $N_{\text{max}}=35$ and tolerance $\varepsilon = 10^{-5}$ with an equispaced training set $\mathcal P_h$ of 51 parameters. We stop the procedure at the iteration for which the condition $\tau_N(\mu) \leq 1$ for every $\mu \in \mathcal P_h$ is verified. In the left plot of Figure \ref{fig:Manzoni}, we show the evolution of $\tau_N(\mu)$ for several iterations of the vanilla-greedy procedure. For this test case, the criterion $\tau_N(\mu) \leq 1$ is met for the $8$-th iteration: at this point the nonlinear estimator $\Delta_N^{\text{nl}}(\mu)$ can be used instead of $\Delta_N^{\text{lin}}(\mu)$.

This is not the case with the numerical approximation of problem (ii), the flow in a sudden-expansion channel depicted in Figure \ref{fig:channel}. To compare the results with test (i), we change the inlet profile to $v_{\text{in}} = [0.96(5-x_2)(x_2 -2.5), 0]^T$ , in order to obtain a characteristic velocity $U=1$ as in \cite{ManzoniEfficientComputationalFramework2014}. The parameter space is given by $\mathcal P = [0.016, 0.064]$, corresponding to $Re \in [39, 156]$, while the discretization is the one proposed in Section \ref{sec:results}. Once again, we employ the vanilla-greedy approach with $N_{\text{max}}=35$ and tolerance $\varepsilon = 10^{-5}$. As we can observe in the right plot of Figure \ref{fig:Manzoni}, the evolution of $\tau_N(\mu)$ has difficulties in reaching the needed threshold. In fact, the condition $\tau_N(\mu) \leq 1$ is never met, and it is still quite far at the $8$-th iteration from reaching such a goal. The vanilla-greedy strategy incurs several problems related to its bifurcating nature such as not convergence of the HF solver and multiple selection of the same parametric instance in the error estimation. For these reasons, in Section \ref{sec:results}, we employ the linear bound $\Delta_N^{\text{lin}}(\mu)$. 
Finally, we remark that, as known in the literature, the linear and nonlinear bounds have comparable behaviors and the accuracy of the method is not affected by this simplification \cite{ManzoniEfficientComputationalFramework2014}.

\end{appendices}
\end{document}